\documentclass[a4, 12pt]{article}
\usepackage{amssymb}
\usepackage{amsmath}
\usepackage{bbm}
\usepackage{empheq}
\usepackage[framemethod=tikz]{mdframed}
\usepackage{lipsum}
\usepackage{tcolorbox}
\usepackage{fourier-orns}
\usepackage{enumitem}
\RequirePackage[amsmath,thmmarks,hyperref,framed]{ntheorem}
\usepackage[T1]{fontenc}
\usepackage{lmodern}
\usepackage{stmaryrd}
\usepackage{cancel}
\usepackage{enumitem}

\usepackage[applemac]{inputenc}

\def \R{\mathbb{R}}

\numberwithin{equation}{section}

\theoremstyle{plain}
\newtheorem{theorem}{Theorem}[section]
\newtheorem{corollary}[theorem]{Corollary}

\newtheorem{lemma}[theorem]{Lemma}
\newtheorem{proposition}[theorem]{Proposition}

\theoremstyle{definition}

\theoremstyle{remark}
\newtheorem{remark}[theorem]{Remark}

\usepackage[english]{babel}
\newmdenv[innerlinewidth=0.5pt, roundcorner=4pt,innerleftmargin=6pt,
innerrightmargin=6pt,innertopmargin=10pt,innerbottommargin=10pt,backgroundcolor=gray!21]{mybox}

\theorembodyfont{\upshape}

\date{}
\title{Ancestral diversity in fragmentation trees}

\author{B\'en\'edicte Haas \thanks{Universit\'e Sorbonne Paris Nord,
    LAGA, CNRS (UMR  7539) 93430 Villetaneuse, France \newline
    \hspace*{0.5cm} E-mail: haas@math.univ-paris13.fr}  \quad \&
  \hspace{0.2cm}  Gr\'egory Miermont \thanks{Ecole Normale Supérieure
    de Lyon, UMPA, CNRS (UMR 5669), 69364 Lyon Cedex 7, France\newline
    \hspace*{0.5cm} Institut
    Universitaire de France\newline
    \hspace*{0.5cm} E-mail: gregory.miermont@ens-lyon.fr}}

\begin{document}

\maketitle

\begin{abstract}
In a deterministic or random tree, a notion of ancestral diversity can
be defined as follows. Sample independently $n$ groups of $k$
leaves and count the number $N_n(k)$ of distinct most recent common ancestors
of each of the groups. As $n$ becomes large, the asymptotic behavior of
$N_n(k)$ depends of course on the structure of the tree. Motivated by
the study of the edge density in the Brownian
co-graphon, Chapuy recently considered this problem
in the case where $k=2$ and where the tree is the Brownian continuum
random tree. We vastly extend this framework by considering general
values of $k$ and general
fragmentation trees, which include some
prominent examples such as stable L\'evy trees and idealized models of
phylogenetic trees. Other natural ancestral statistics are also considered.
For a given tree model, we identify a phase
transition-like phenomenon, with 
different asymptotic regimes for $N_k(n)$, depending on the position
of $k$ relative to a model-dependent critical value. 
\end{abstract}

\section{Introduction and main results}

This paper is concerned with certain statistics associated with random
tree-like structures. To fix the ideas, let $(T,d)$ be an $\R$-tree, rooted at a distinguished point
$\rho\in T$, and let $\mu$ be a nonatomic probability measure on $(T,d)$, which
is supported on the set of leaves of $T$, that is, on points $x\in T$
such that $T\setminus \{x\}$ is connected.
Fix two integers $k\geq 2$ and
$n\in \mathbb{N}=\{1,2,\ldots\}$. Let $(x_{i,j},1\leq
i\leq n,1\leq j\leq k)$ be a $n\times k$ array of independent random
variables with law $\mu$.
For every $i\in \{1,2,\ldots,n\}$, we let $a_i(k)$ be the
most recent common ancestor of $x_{i,1},\ldots,x_{i,k}$, which is the
unique point $y$ such that the segment in $T$ from $\rho$ to $y$ is
the intersection of the segments from $\rho$ to $x_{i,j},1\leq j\leq k$. We are
concerned with the behavior of the \emph{ancestor-counting} random variable
$$N_n(k)=\mathrm{Card}\{a_i(k):1\leq i\leq n\}$$
as $n$ converges to infinity. 

It might be the case that this random variable is of some interest in
phylogenetics, as it could be used as a measure of genetic diversity
in a given population. However, our motivation for studying this
problem does not come from mathematical biology, but rather from a paper by
Chapuy \cite{chapuy25}, who was interested in analyzing the
moments of the edge-density of the Brownian co-graphon introduced in
\cite{BaBoFeGeMaPi22}. He considers the case where the tree $T$ is the
Brownian continuum random tree and where $k=2$. Combinatorially, the
problem is equivalent to the following: let $B_n$ be a uniform random
rooted plane binary tree with $2n$ leaves labeled by $1,2,\ldots,2n$,
and let $b_i$ be the most recent common ancestor to the leaves labeled
$i$ and $i+n$. Then $N'_n=\mathrm{Card}\{b_i:1\leq i\leq n\}$ has
same distribution as the random variable $N_n(2)$. Chapuy shows that
$N_n(2)/(\sqrt{n}\log n)$ converges in $L^2$ to $1/\sqrt{2\pi}$, by using 
a second moment method, and gives an explanation to the logarithmic
factor by an argument that has some ``analytic number theoretic''
flavor. He also asks if the limiting behavior of $N_n(2)$ could be
derived by other means that would involve natural processes related to
the structure of the tree.

Our paper provides such a derivation, 
and also extends Chapuy's result in various ways. First, we allow the
tree structure to belong to a larger family of
trees with a certain fragmentation property, first considered in
\cite{HaMi04}, as we will recall in Section
\ref{sec:fragmentation-trees-2}.
Second, we allow for $k$-tuples of sampled vertices with 
arbitrary $k\geq 2$. Finally, we  also  consider other natural
ancestor-counting statistics such
as the number $N_{n,r}(k)$ of ancestors $a$ with multiplicity $r$,
that is, such that $\mathrm{Card}\{i\in 
\{1,2,\ldots,n\}:a_i(k)=a\}=r$.  

Our approach, valid for any rooted $\mathbb R$-tree equipped with a
nonatomic probability measure on its set of leaves, is to represent
the random variable $N_n(k)$ as the number 
of distinct boxes in an urn process, and to apply classical results of
Karlin \cite{karlin67,GPY06}, relating this number to the asymptotic
behavior of the number of urns exceeding a given size $x$ as $x\downarrow
0$. As it turns out, in our general context of fragmentation trees,
these urn-counting random variables arise as a
particular instance of {\em large dislocations in a self-similar
  fragmentation}, as considered in a work by Quan Shi \cite{qshi15},
which extended earlier results by Bertoin and Mart\'inez
\cite{BeMa05}. Interestingly, \cite{qshi15} showed that a  phase 
transition-type phenomenon occurs depending on the fragmentation
mechanism. 
The general phenomenon that
we observe is that there exists a critical value $k_0$, which depends
on the law of the tree, and is not
necessarily an integer, such that if $k>k_0$, then the properly
renormalized urn count converges to a nondeterministic limit, while if $k<  k_0$,
then it admits a deterministic scaling limit. 

However, the results of \cite{qshi15}, do not apply in a
direct way to our setting for two reasons. First, they provide
convergence in $L^2$ for the rescaled urn counts, while Karlin's
result requires almost sure convergence in order to transfer these
results to the ancestor-counting random variables $N_n(k)$. For this
reason, we need to quantify the speed of convergence in $L^2$, which
was not addressed in \cite{qshi15}. Second, they
do not encompass the critical case $k=k_0$ (see \cite[Remark 
2.7]{qshi15}), which is precisely the situation of Chapuy's
result. 
Hence we resume, in a sense, where 
\cite{qshi15} stopped, and introduce new techniques to deal with these
issues. In particular, in the critical case $k=k_0$,
(which requires that $k_0$ be an integer) we will see
that logarithmic
corrections arise, but the limit is still deterministic, as in
Chapuy's result.

Before presenting
the general method, let us discuss in more details the situation in
the particular case of the Brownian continuum random tree.

\subsection{The Brownian tree case}\label{sec:brownian-tree-case}

In \cite{chapuy25}, Chapuy considered the number of distinct
ancestors of a sample of $n$ pairs of leaves in the Brownian CRT. In
this work, we show that his result can be generalized to the case
where one picks $k$-tuples of leaves at once, for some fixed integer $k\geq 2$.
We observe that the situation is very different if $k=2$ or $k\geq 3$,
in the sense that the limits in the latter case are random, while they
are deterministic in the former case considered in \cite{chapuy25}. 

Let $\mathbf{e}$ be a normalized Brownian excursion. For every $t\geq
0$, we let $L_1(t)\geq L_2(t)\geq \ldots\geq 0$ be the ranked sequence of
Lebesgue measures of the connected components of the open set
$\{x\in [0,1]:\mathbf{e}_x>t\}$. For every $k\geq 3$, we define the
random variable
$$X_k=2\sqrt{2k}\, \cdot \int_0^\infty \sum_{i\geq 1}L_i(t)^{(k-1)/2}\mathrm dt\, .$$

\begin{theorem}
For the Brownian CRT, it holds that, almost surely and in $L^2$,
$$
\lim_{n\to\infty}\frac{N_n(2)}{\sqrt{n}\, \log (n)} =\frac{1}{\sqrt{2 \pi}} \ ;
\qquad \lim_{n\to\infty} \frac{N_n(k)}{\sqrt{n}}  = X_k \quad \text{
  for } k \geq 3\ .
$$
Moreover,  for every $r\geq 1$, we have the following almost sure limits
$$
\lim_{n\to\infty}\frac{N_{n,r}(2)}{\sqrt{n}\, \log(n)}  =\frac{
  \Gamma(r-1/2)}{2\sqrt 2 r! \pi} 
 \ ; \qquad \lim_{n\to\infty} \frac{N_{n,r}(k)}{\sqrt{n}}  =\frac{
  \Gamma(r-1/2)}{2 r! \sqrt{\pi}} \cdot X_k \quad \text{ for } k \geq 3.
$$
\end{theorem}

\smallskip

{\remark We note that $X_3$ is $2 \sqrt{6}$ times the area under the  standard
Brownian excursion of length 1, which has a well-studied law sometimes
called the Airy area distribution, see \cite{janson07}. In particular,
the moments of $X_3$ solve explicit quadratic recurrence equations,
see formulas (4--9) therein.  This holds in fact for
every value of $k\geq 3$, by \cite[Corollary 2.2]{bertoin12area}, see
the discussion before Theorem \ref{sec:thm-super} below. In
particular, the first moment admits the
expression
$$\mathbb{E}[X_k]=\frac{\sqrt{k}\,
 \cdot  \Gamma(k/2-1)}{\Gamma((k-1)/2)}\, .$$}

\subsection{Fragmentation trees}\label{sec:fragmentation-trees-2}

A self-similar fragmentation process \cite{berthfrag01,bertsfrag02} describes the
evolution of a system of massive objects which are subject to a random
splitting as time evolves. Informally, the system starts from a single
object of mass $1$, and at any given time, an object of size $x$ is
dislocated into sub-objects of sizes $x\mathbf{s}=(xs_1,xs_2,\ldots)$ at a rate
$x^{\alpha}\nu(\mathrm d\mathbf{s})$, where $\alpha$ is a real number
and $\nu$ is a {\em dislocation measure}. This means that 
$\nu$ is a $\sigma$-finite measure on the set
$$\mathcal{S}=\left\{\mathbf{s}=(s_1,s_2,\ldots):s_1\geq s_2\geq
\ldots \geq 0\, ,\sum_{i\geq 1}s_i\leq 1\right\}\, ,$$
which satisfies $\nu(\{(1,0,0,\ldots)\})=0$, as well as the integrability
condition
$\int_{\mathcal{S}}(1-s_1)\nu(\mathrm d\mathbf{s})<\infty$. To
avoid trivialities, we will always assume that $\nu(\mathcal{S})>0$,
and we will also make the simplifying assumption that
$\nu(\{\mathbf{s}\in \mathcal{S}:\sum_{i\geq 1} s_i<1\})=0$,
so that the total mass is preserved at each dislocation event. In \cite{berthfrag01,bertsfrag02}, Bertoin proved that for every such 
$(\alpha,\nu)$, there exists a process
$(F^{(\alpha,\nu)}_i(t),i\geq 1)_{t\geq 0}$ with values in $\mathcal{S}$,
which evolves according to the Markovian dynamics heuristically
described above. 

When $\alpha<0$, it was shown in
\cite{HaMi04} that the process can be described in terms of
a random compact measured rooted $\R$-tree $(T^{(\alpha,\nu)},d,\rho,\mu)$, in the
sense that the process $F^{(\alpha,\nu)}$ has the same distribution as
$$\mu(T^{(\alpha,\nu)}_i(t)),i\geq 1,\quad t\geq 0\, ,$$
where, for every $t\geq 0$, the sets $T^{(\alpha,\nu)}_i(t),i\geq 1$ are the
connected components of $\{x\in T^{(\alpha,\nu)}:d(x,\rho)>t\}$, indexed by
decreasing order of their $\mu$-measures. In particular, we can
associate with the tree $(T^{(\alpha,\nu)},d,\rho,\mu)$
the ancestor-counting random variables $N^{(\nu)}_n(k),N^{(\nu)}_{n,r}(k)$ of the
introduction. As the notation suggests, the laws of these random
variables are actually independent of $\alpha$, as will be discussed in
Section~\ref{sec:part-repr-gene}. 

This framework encompasses the case of the Brownian CRT, which is
obtained for 
\begin{equation}
\label{nu_brownien}
\alpha=-\frac{1}{2}\, ,\qquad
\int_{\mathcal{S}}\nu(\mathrm{d}\mathbf{s})f(\mathbf{s})=\sqrt{\frac{2}{\pi}} \cdot \int_{1/2}^1\frac{\mathrm{d}
x}{(x(1-x))^{3/2}} f(x,1-x,0,0,\ldots)\, .
\end{equation}
Many other classes of random continuum trees can be obtained in this
way, see the discussion of Section~\ref{sec:examples}. Notably, we
emphasize that, like the Brownian CRT, these models of  fragmentation
trees appear as scaling limits of many natural models of discrete
trees \cite{HaMi12,HaMiPiWi08}. See also the recent monograph \cite{BCR-SSMT} for
a generalization to the framework of self-similar Markov trees. 

\subsection{Main results}\label{sec:main-results}

We fix a dislocation measure $\nu$, and let $\gamma\in
[0,1)$. Let us make the following assumption, which will be key in all the
results discussed in this paper. 
\begin{equation}
  \label{eq:hg}
\mbox{There exists } c_\nu\in (0,\infty) \mbox{ such that } \nu(s_1
\leq 1-x) \underset{x \downarrow 0} \sim c_\nu x^{-\gamma}\, .
  \tag{$\mathbf
{H_{\gamma}}$}  
\end{equation}

\subsubsection{Supercritical case}

We let $\nu$ be a dislocation measure satisfying \eqref{eq:hg}, and
first consider an integer $k$ such that $k\gamma>1$. We call this situation
the {\em supercritical case}. Note that this requires in particular
that $\gamma>0$. We let $(T^{(1-k\gamma,\nu)},d,\rho,\mu)$ be the
self-similar fragmentation tree with index $\alpha=1-k\gamma$ and
dislocation measure $\nu$, as discussed in Section
\ref{sec:fragmentation-trees-2}, 
and define the following random variable
\begin{equation}
  \label{eq:23}
  A^{(\nu)}_k=c_\nu\int_{T^{(1-k\gamma,\nu)}}d(\rho,u)\mu(\mathrm
  du) \, .
\end{equation}
If, as discussed in Section
  \ref{sec:fragmentation-trees-2}, we define a fragmentation process by letting 
$F^{(1-k\gamma,\nu)}(t)$ be the decreasing sequence of $\mu$-measures
of the connected components of $\{x\in
T^{(1-k\gamma,\nu)}:d(\rho,x)>t\}$, then we have the alternative
formula
$$A^{(\nu)}_k=c_\nu\int_0^\infty\sum_{i\geq
  1}F^{(1-k\gamma,\nu)}_i(t)\, \mathrm dt\, ,$$
which, up to the factor $c_\nu$, is called the {\em area} of the
fragmentation process in \cite{bertoin12area}. 
When $\nu$ is binary, meaning that $\nu(\{\mathbf{s}\in
\mathcal{S}:s_3>0\})=0$,  
the moments of this random variable satisfy certain explicit
quadratic recursive formulas, as shown in \cite[Corollary 2.2]{bertoin12area}. 
Finally, we note that this variable is ``homogeneous'', in the sense that
$A_k^{(\lambda \nu)}$ has same distribution as $A_k^{(\nu)}$ for every
$\lambda>0$. This comes from the fact that $F^{(1-k\gamma,\lambda
  \nu)}(\cdot/\lambda)$ has the same distribution as
$F^{(1-k\gamma,\nu)}$. 

\begin{theorem}[Supercritical case, $k\gamma>1$]\label{sec:thm-super}
  Assume \eqref{eq:hg} and let $k$ be such that $k\gamma>1$. Then the following
  limit holds almost surely and in $L^2$:
$$
\lim_{n\to\infty}\frac{N^{(\nu)}_n(k)}{ n^{\gamma}} =  \Gamma(1-\gamma)
k^{\gamma}  \cdot A_k^{(\nu)}\, .$$
Moreover, we have, for every $r\geq 1$, almost surely: 
$$\lim_{n\to\infty}\frac{N^{(\nu)}_{n,r}(k)}{n^{\gamma}}  =\frac{\gamma
  \Gamma(r-\gamma) k^{\gamma} }{r!} \cdot A_k^{(\nu)}\, .
$$
\end{theorem}

\medskip

{\remark We believe that the last stated convergence also holds in
  $L^2$, but we haven't checked the details. A similar remark applies
  to the forthcoming Theorem \ref{sec:thm-sub-gammapos}. }

{\remark 
  The expectation of $ A_k^{(\nu)}$ is equal to $c_\nu/\phi(k\gamma-1)$, see Section \ref{sec:tagged-fragment}. Therefore, under \eqref{eq:hg}, the random variable appearing as the limit of $n^{-\gamma} \cdot N^{(\nu)}_n(k)$ in the supercritical regime has an expectation which converges as $k\rightarrow \infty$: 
$$\mathbb E\big[\Gamma(1-\gamma)  k^{\gamma} \cdot A^{(\nu)}_k\big] \underset{k\rightarrow \infty }{\longrightarrow} \gamma^{-\gamma}.$$}

\subsubsection{Subcritical and critical cases}\label{sec:subcritical}

Still working under \eqref{eq:hg}, we now assume that $k\geq 2$ is
such that $k\gamma\leq 1$. We call this situation the {\em subcritical
case} when $k\gamma<1$, and the {\em critical case} when
$k\gamma=1$.

We consider the following assumption. 
\begin{equation}
  \label{eq:exp}
\mbox{ There exists } \eta\in (0,1) \mbox{ such 
that }  \int_{\mathcal S}\sum_{i\geq 2}s_i^{1-\eta} 
  \, \nu(\mathrm d \mathbf s)<\infty\, . 
\tag{\textbf{Exp}}  
\end{equation}
Note that this is automatically verified
if $\nu(s_{m+1}>0)=0$ for some $m\geq 2$. 
We also consider one last assumption that will be useful in the
case
$\gamma=0$. 
\begin{equation}
  \label{eq:gamma0}
  \mbox{The measure }\sum_{i\geq 1}\nu(s_i\in \mathrm d x)\mbox{ is 
    absolutely continuous.}
  \tag{\textbf{Dens}}  
\end{equation}
Finally, 
we define 
\begin{equation}
  \label{eq:24}
  C^{\mathrm{sub}}_\nu(k)= \frac{\int_{\mathcal S}\nu(\mathrm d\mathbf{s})\big(1-\sum_{i\geq 1}
  s_i^k\big)^{\frac{1}{k}}}{\int_{\mathcal{S}}\sum_{i\geq
    1}s_i|\log(s_i)|\nu(\mathrm d\mathbf{s})}\quad \mbox{ and }\quad
C^{\mathrm{cr}}_\nu(k)=\frac{c_\nu\, k^{\frac{1}{k}-1}}{\int_{\mathcal{S}}\sum_{i\geq
      1}s_i|\log(s_i)|\nu(\mathrm d\mathbf{s})}\, .
\end{equation}
Under \eqref{eq:exp}, $~C^{\mathrm{cr}}_\nu(k)$ is finite, positive and
a homogeneous functions of $\nu$, and so does $C^{\mathrm{sub}}_\nu(k)$ when $k\gamma<1$.

\vspace{0.1cm}

\begin{theorem}[Subcritical and critical cases, $k\gamma\leq 1$]
  \label{sec:thm-sub-gammapos}
  Let us assume that \eqref{eq:hg} holds for some $\gamma\in [0,1)$,
  that \eqref{eq:exp} also holds, and let $k\geq 2$. If $\gamma=0$, we also
  assume that \eqref{eq:gamma0} holds. 

  In the subcritical case $k\gamma<1$, it holds that, almost
  surely and in $L^2$, 
$$
\lim_{n\to\infty} \frac{N^{(\nu)}_n(k)}{n^{\frac{1}{k}}} = \Gamma\left(1-\frac{1}{k}\right)C^{\mathrm{sub}}_\nu(k)$$
In the critical case $k\gamma=1$, it holds that, almost surely and in
$L^2$, 
$$
\lim_{n\to\infty}\frac{N^{(\nu)}_n(k)}{n^{\frac{1}{k}} \log(n)} =
\Gamma\left(1-\frac{1}{k}\right)  C^{\mathrm{cr}}_\nu(k)\, .
$$
Moreover, for every $r\geq 1$, we have the following almost sure
limits, respectively when $k\gamma<1$ and $k\gamma=1$: 
$$\lim_{n\to\infty}
\frac{N^{(\nu)}_{n,r}(k)}{n^{\frac{1}{k}}}=
\frac{\Gamma(r-\frac{1}{k})}{k\, r!} C^{\mathrm{sub}}_\nu(k)
\qquad \mbox{ and }\qquad  \lim_{n\to\infty}\frac{N^{(\nu)}_{n,r}(k)}{n^{\frac{1}{k}} \log(n)}=
\frac{\Gamma(r-\frac{1}{k})}{k\, r!} C^{\mathrm{cr}}_\nu (k)
  \, .
  $$
\end{theorem}

The assumptions \eqref{eq:exp} and \eqref{eq:gamma0} made in this
statement are certainly not optimal, but hold in all 
the examples discussed in the paper. For instance, we could weaken
\eqref{eq:gamma0} a little bit, by assuming  that
the measure, or some multiplicative convolution thereof, has a
non-trivial absolutely continuous part. 

\medskip

{\remark[Ancestors multiplicities]{
Our approach also yields immediately the following result, valid for
any rooted measured $\R$-tree $(T,d,\rho,\mu)$}. Almost surely and in $L^1$,
$$
\sum_{b\in \mathrm{Br}(T)} \left|\frac{D_{b,n}(k)}{n}-P(k,T,b)
\right|\underset{n\to\infty}{\longrightarrow} 0\, ,
$$
where $\mathrm{Br}(T)$ is the set of branchpoints of $T$, that is, of
points $b$ such that $T\setminus \{b\}$ is not connected, $D_{b,n}(k)$ is the number of $k$-samples, amongst the $n$
first, that have the branch point $b$ as most recent common ancestor,
and $P(k,T,b)$ is defined in \eqref{eq:3}. In particular, almost surely,
$$
\frac{1}{n}\, \max_{b\in \mathrm{Br}(T)}D_{b,n}(k)\underset{n\to\infty}{\longrightarrow} \max_{b\in \mathrm{Br}(T)} P(k,T,b) \, .
$$} 

\subsection{Organisation of the paper}

Section \ref{sec:paintbox} reformulates our problem in terms of urn models and recalls Karlin's classical result on counting occupied urns. We also review there some framework on self-similar fragmentation processes and trees, as well as  elements of renewal theory for subordinators, in relation with the tagged fragment process. In order to apply Karlin's result, we need to study a notion of ``large'' dislocations in fragmentation processes, as already considered by Quan Shi's in \cite{qshi15}. To this end, in Section \ref{sec:sizes-small-boxes}, we lay out the first steps of our approach, based on the 
key renewal theorem for subordinators and a first concentration inequality, following a similar line to \cite{qshi15}. Theorem \ref{sec:thm-super} on the supercritical case then follows rather easily and its proof is included in Section \ref{sec:sizes-small-boxes}. The subcritical and critical cases of Theorem \ref{sec:thm-sub-gammapos} are more involved and studied in Section \ref{sec:subcritical-case}. They rely on a second concentration inequality, which is more difficult to establish and requires finer renewal estimates. Finally, Section \ref{sec:examples} is devoted to several examples of applications, notably the stable L\'evy trees of Duquesne, Le Gall and Le Jan (including the Brownian CRT), and two one-parameter families of theoretical models for phylogenetic trees: Ford's model and Aldous's beta-splitting model.

\medskip

\noindent{\bf Acknowledgements. }Part of this research was conducted
while GM was holding a visiting professor position at the Research
Institute for Mathematical Sciences, an International Joint
Usage/Research Center located in Kyoto University.

\section{Preliminaries}
\label{sec:paintbox}

This section lays the basics of our approach, which consists in
reformulating our problem in terms of classical urn schemes, and to
express this urn scheme in terms of the appropriate statistics of the
fragmentation processes. 

\subsection{Reformulation as an urn-counting
  problem}\label{sec:urn-schemes} 

\subsubsection{Classical urn schemes. }
Let $\mathbf{p}=(p_1,p_2,\ldots)$ be a nonincreasing
sequence of nonnegative numbers with sum $1$. Let
$\xi_1,\xi_2,\ldots$ be an i.i.d.\ sequence of random variables with
law $\mathbf{p}$: we imagine that a ball labeled
$i$ falls into an urn with label $j$ with probability $p_j$. 

Let
$$N_n=\mathrm{Card}\left(\{\xi_1,\ldots,\xi_n\}\right)$$
be the number of nonempty urns after $n$ draws, and
$$N_{n,r}=\mathrm{Card}\left(\left\{j\geq 1:\sum_{i=1}^n\mathbbm{1}_{\{\xi_i=j\}}=r \right\}\right)$$
be the number of urns containing exactly $r$ balls after $n$ draws. We
call these the {\em urn-counting} random variables associated with
$\mathbf{p}$. 

A famous work of Karlin \cite{karlin67} shows that the
asymptotic behavior of the random variables $N_n,N_{n,r}$ is
intimately linked with the decrease rate of $p_i,i\geq 1$, expressed
in terms of the {\em urn distribution function} 
\begin{equation}
  \label{eq:2}
  S_x^{\mathbf{p}}=\max\{j\geq 1:p_j\geq
  x\}\, ,\qquad  x>0\, .
\end{equation}
We state an improved form of this result, due to
Gnedin-Pitman-Yor \cite[Theorem 2.1]{GPY06}, that allows the sequence
$\mathbf{p}$ to be itself random, in which case the urn scheme
described above and the random variables $\xi_i,N_n,N_{n,r}$ are all 
defined conditionally on $\mathbf{p}$. 

\begin{theorem}
\label{thm:Karlin}
Let $\mathbf{p}$ be a random nonincreasing sequence with sum $1$.
Assume that there exist a real number $\rho\in (0,1)$, a function $\ell$
that is slowly varying at $\infty$, and a nonnegative random variable $L$, such
that $\lim_{x\downarrow
  0}\frac{x^{\rho}}{\ell(1/x)}S^{\mathbf{p}}_x=L$ almost surely. 
Then it holds that
$$\lim_{n\to\infty} \frac{N_n}{n^{\rho}\ell(n)}=\Gamma(1-\rho)L\,
,\qquad \lim_{n\to\infty}
\frac{N_{n,r}}{n^{\rho}\ell(n)}=\frac{\rho\Gamma(r-\rho)}{r!}L\, ,$$
almost surely.
\end{theorem}

\subsubsection{Reformulation of the ancestor-counting random variables}\label{sec:reform-ancest-count}

Let us now reformulate the ancestor-counting random variable $N_n(k)$ of
a tree in terms of an urn scheme. Let $(T,\rho, d,\mu)$ be a
compact, rooted, measured $\mathbb{R}$-tree, with $\mu$ a nonatomic
probability measure that only charges the set of leaves of $
T$. We let $\mathrm{Br}(T)$ be the set of branchpoints of $T$, that
is, the set of points $b$ such that $T\setminus \{b\}$ has at least
two connected components not containing $\rho$. This set is at most
countable, and for $b\in\mathrm{Br}(T)$, we let $T_b$ be the union of
all connected components $T_{b,i},i\geq 1$ of $T\setminus \{b\}$ not
containing $\rho$, where the latter are labelled by nonincreasing
order of $\mu$-measure.

Let $X_1,\ldots,X_k$ be i.i.d.\ with distribution $\mu$. 
Clearly, an element $b\in \mathrm{Br}(T)$
is the common ancestor of $X_j, 1\leq j\leq k$ if and
only if these $k$ leaves are all in $T_b$, but not all in a
common subtree $T_{b,i}$ for some $i\geq 1$.
This event occurs with probability $\mu(T_b)^k-\sum_{i\geq
  1}\mu( T_{b,i})^k$, and since $\mu$ is nonatomic and supported on the
leaves of $T$, these probabilities sum to $1$ as $b$ describes
$\mathrm{Br}(T)$.
Therefore, if we let $P_1\geq P_2\geq
\ldots\geq 0$ be the nonincreasing rearrangement of the family 
\begin{equation}
  \label{eq:3}
  P(k,T,b)=\mu( T_b)^k-\sum_{i\geq
    1}\mu( T_{b,i})^k\, ,\qquad b\in \mathrm{Br}(T)\, ,
\end{equation}
then the ancestor-counting random variables
$N_n(k),N_{n,r}(k)$ are nothing but the urn-counting random variables
associated with the sequence $\mathbf{P}=(P_1,P_2,\ldots)$. 

When the tree $(T,d,\rho,\mu)$ is a fragmentation tree, we may
re-express the  associated urn count distribution process in terms of the
associated fragmentation process, as we will now see.

\subsection{Basic tools of fragmentation 
  processes}\label{sec:basic-tools-fragm} 

Let $\nu$ be a dislocation measure, and $\alpha$ a real number,
associated with a self-similar fragmentation $F^{(\alpha,\nu)}$.

\subsubsection{Partition representation and genealogy}\label{sec:part-repr-gene}

According to the discussion of Section
\ref{sec:fragmentation-trees-2}, one can
view a fragmentation process as an $\mathcal{S}$-valued process
recording the masses of the objects
present at time $t$. However, with this point of view, the natural
genealogical structure of the process is lost. A similar situation is
classically encountered in the study of branching processes, where one
can focus only on the evolution of the total population size, or
consider the genealogical tree of the population as well.

The key idea of Bertoin
\cite{berthfrag01,bertsfrag02} is to represent a fragmentation process as
a process $(\Pi(t),t\geq 0)$ with values in the set of
partitions of $\mathbb{N}$, which is nondecreasing in the sense that
$\Pi(t)$ is finer than $\Pi(s)$ for every $t\geq s\geq 0$, and 
whose law is exchangeable, that is,
invariant under the action of the permutation group of
$\mathbb{N}$. He showed that  for every $(\alpha,\nu)$ as
above, there is a unique (in law) such process
$\Pi=(\Pi^{(\alpha,\nu)}(t),t\geq 0)$ such that almost surely, for every $t\geq 0$, every block $B$
of the partition $\Pi^{(\alpha,\nu)}(t)$ admits an asymptotic frequency
$$|B|=\lim_{n\to\infty}\frac{\mathrm{Card}(B\cap\{1,2,\ldots,n\})}{n}\,
,$$ and such that 
the process $(F^{(\alpha,\nu)}_i(t),i\geq 1)_{t\geq 0}$ of these
asymptotic frequencies, ranked in nonincreasing order, obeys the
Markovian dynamics heuristically described in Section
\ref{sec:fragmentation-trees-2}, that is, every object of size $x$ 
dislocates into sub-objects of sizes $x\mathbf{s}$ at infinitesimal
rate $x^\alpha\nu(\mathrm d\mathbf{s})$. 

Moreover, if we let $\Pi_{(i)}^{(\alpha,\nu)}(t)$
denote the block of $\Pi^{(\alpha,\nu)}(t)$ containing the integer
$i$, then we may couple the processes $\Pi^{(\alpha,\nu)}$ together,
for a fixed choice of $\nu$, in such a way that, for every $\alpha\in
\R$, $t\geq 0$ and $i\in \mathbb{N}$, 
\begin{equation}
  \label{eq:29}
  \Pi_{(i)}^{(\alpha,\nu)}(t)=\Pi_{(i)}^{(0,\nu)}(\tau_{(i)}^{(\alpha)}(t))\,
  ,
\end{equation}
where
\begin{equation}
  \label{eq:28}
  \tau^{(\alpha)}_{(i)}(t)=\inf\left\{s\geq 0:
    \int_0^s|\Pi^{(0,\nu)}_{(i)}(u)|^{-\alpha}\mathrm du>t\right\}\, .
\end{equation}

The important feature of this representation is that it is now
possible to associate a genealogy to the process
$(\Pi^{(\alpha,\nu)}(t),t\geq 0)$. In fact, when $\alpha<0$, then
\cite{HaMi04} showed that there exists a unique (in law) random rooted
and 
measured $\R$-tree $(T^{(\alpha,\nu)},d,\rho,\mu)$ such that, if
$x_i,i\geq 1$ is an independent sample of $\mu$-distributed random
points, then the partition-valued process $(\Pi(t),t\geq 0)$ defined
by the property that $i$ and $j$ are in the same block of $\Pi(t)$ if
and only $x_i,x_j$ belong to the same connected component of $\{x\in
T^{(\alpha,\nu)}:d(x,\rho)>t\}$, has the same distribution as
$\Pi^{(\alpha,\nu)}$. For this reason, we may and will actually assume
that $\Pi=\Pi^{(\alpha,\nu)}$. 

In this representation, the 
branchpoints $b$ of $T^{(\alpha,\nu)}$
correspond to the dislocation events in the process
$\Pi^{(\alpha,\nu)}$, that is, the pairs $(B,t)$ such that $B$ is a block of
$\Pi^{(\alpha,\nu)}(t-)$ (which is the coarsest partition that
is finer than $\Pi^{(\alpha,\nu)}(s)$ for every $s<t$), 
but not a block of $\Pi^{(\alpha,\nu)}(t)$.  Moreover, with this
correspondence, and using the notation around
\eqref{eq:3},
$$\mu(T_b)=|B|\, ,\qquad \mu(T_{b,i})=|B_i|\, ,$$
where $B_1,B_2,\ldots$ are the blocks of $B\cap \Pi^{(\alpha,\nu)}(t)$, arranged by
decreasing order of asymptotic frequency. Finally, because of the
correpondence \eqref{eq:29}, we may and will assume that
$\Pi^{(\alpha,\nu)}$ is associated with a homogeoneous fragmentation
process $\Pi^{(0,\nu)}$. 

\begin{proposition}
  \label{sec:part-repr-gene-1}
  For every $\alpha<0$, the urn sizes \eqref{eq:3} associated with the self-similar
  fragmentation tree $T^{(\alpha,\nu)}$ are equal to 
  the decreasing rearrangement of the family
  $$|\Pi_i^{(0,\nu)}(t-)|^k-\sum_{j\geq
    1}|\Pi_{i,j}^{(0,\nu)}(t)|^k\, ,\quad t\geq 0,i\in \mathbb N\, ,$$
  where $\Pi^{(0,\nu)}_{i,j}(t),j\geq 1$ are the blocks of
  $\Pi^{(0,\nu)}(t)$ that are contained in $\Pi^{(0,\nu)}_i(t-)$. 
\end{proposition}

Note that the resulting law does not depend on $\alpha$. This is due
to the fact that the mass of subtrees does not depend on the tree
metric, but only on the genealogical structure. For this reason, we
now work exclusively with homogeneous fragmentations
$\Pi^{(\nu)}=\Pi^{(0,\nu)}$ and $F^{(\nu)}=F^{(0,\nu)}$.

\subsubsection{The tagged fragment}\label{sec:tagged-fragment}

      An auxillary process of crucial importance is 
the {\em tagged fragment} process $(F_*(t)=|\Pi_{(1)}^{(\nu)}(t)|$, $t\geq 0)$, which can be 
seen as the size at time $t$ of the object containing 
a point marked uniformly at random according to the total mass measure. 
  It satisfies the following many-to-one formula: for every 
  measurable $f:\R\to 
  \R_+$ and $t\geq 0$, 
  $$\mathbb{E}\left[\sum_{i\geq 1}F^{(\nu)}_i(t)f(F^{(\nu)}_i(t))\right]=\mathbb{E}[f(F_*(t))]\, .$$
It can also be written as  
$(F_*(t)=e^{-\xi_*(t)},t\geq 0)$, where $\xi_*$ is a subordinator with Laplace 
exponent 
\begin{eqnarray*}
\phi(q)&=&-\log\mathbb{E}[\exp(-q\xi_*(1))] \\
&=&\int_{\mathcal{S}}\Big(1-\sum_{i\geq 
    1}s_i^{q+1}\Big)\nu(\mathrm d \mathbf{s}) \\
    &=&\int_{(0,\infty)}(1-e^{-qx})\Xi(\mathrm d 
  x)\, \, ,\qquad q\geq 0\, , 
 \end{eqnarray*} 
  where $\Xi(\mathrm dx)=\sum_{i\geq 1}e^{-x} \nu(-\log(s_i)\in \mathrm dx)$ 
is the Lévy measure, see \cite{berthfrag01}. 
Our working assumptions admit natural interpretations in terms of
these objects. 

\begin{lemma}
  \label{sec:tagged-fragment-1}
  \begin{itemize}
   \setlength{\itemsep}{0.1pt}
  \item 
Assumption \eqref{eq:hg} is equivalent to
\begin{equation}\label{eq:13}
\phi(q) \underset{q \rightarrow \infty} \sim \Gamma(1-\gamma)c_\nu
q^{\gamma}
\end{equation}
\item 
Assumption \eqref{eq:exp} is equivalent
 to $\phi$ admitting an 
analytic continuation in $(-\eta,\infty)$ for some $\eta>0$. 
\item 
Assumption \eqref{eq:gamma0} is
equivalent to $\Xi(\mathrm dx)$ being absolutely continuous.
\end{itemize}
\end{lemma}

\noindent{\bf Proof. }
For the first point, we note that for all $q\geq 0$
$$
\int_{\mathcal S} \sum_{i\geq 2} s_i^{q+1} \nu(\mathrm d \mathbf s) ~\leq ~2^{-q}\int_{\mathcal S} \sum_{i\geq 2} s_i \nu(\mathrm d \mathbf s) = 2^{-q}\int_{\mathcal S} \left(1-s_1 \right) \nu(\mathrm d \mathbf s),
$$
leading to
$
\phi(q)=\int_{\mathcal S}(1-s_1^{q+1}) \nu(\mathrm d \mathbf s)+O(2^{-q}). 
$
We then conclude with an integration by parts. 

 For the second point, we simply observe that 
   \begin{equation}
     \label{eq:27}
        \int_0^\infty(e^{\eta x}-1)\Xi(\mathrm dx)=\int_{\mathcal{S}}
   \Big(\sum_{i\geq 1}s_i^{1-\eta}-1\Big)\nu(\mathrm 
   d\mathbf{s}) =\int_{\mathcal{S}}
   \sum_{i\geq 2}s_i^{1-\eta}\nu(\mathrm 
   d\mathbf{s})+\int_{\mathcal{S}}
   (1-s_1^{1-\eta})\nu(\mathrm 
   d\mathbf{s})  \, , 
 \end{equation}
 where the last integral is always finite because of the assumption
 that $\int_{\mathcal{S}}(1-s_1)\nu(\mathrm
 d\mathbf{s})<\infty$.

 The third point is immediate. 
 \hfill$\square$

\subsubsection{Potential and resolvent measures}\label{sec:potent-resolv-meas}

A key element of our analysis is the renewal theorem for the potential measure of the
subordinator $\xi_*$, which we now introduce. 
For $\lambda\geq 0$, we let $U_\lambda$ be the 
$\sigma$-finite measure on $\R_+$ defined by 
\begin{equation*}
\int_{\mathbb R_+}
f(y)U_\lambda(\mathrm{d}y)=\mathbb{E}\left[\int_0^\infty e^{-\lambda t}
  f(\xi_*(t))\mathrm{d}t\right]
\end{equation*}  
for every measurable $f:\R_+\to \R_+$. In particular, it is 
characterized by its Laplace transform 
\begin{equation}
  \label{eq:32}
  \mathcal{L}_\lambda(q)= \int_{\mathbb R_+} e^{-qy}U_\lambda(\mathrm{d}y)=\frac{1}{\lambda
  +\phi(q)}\, .
\end{equation}
For $\lambda>0$, $U_\lambda$ has mass $1/\lambda$ and is called the
resolvent measure. Note that $U_1$ is a probability distribution with mean
\begin{equation}
  \label{eq:25}
  \int_{\mathbb R_+} xU_1(\mathrm dx)=\phi'(0+)=\int_{\mathcal{S}}\sum_{i\geq
  1}s_i|\log(s_i)|\nu(\mathrm d\mathbf{s})\, ,
\end{equation}
which is the denominator of the constants \eqref{eq:24}.

On the other hand, the infinite measure $U_0=U$ is called the
potential measure. It has the property that $a\mapsto U([0,a])$ is
a subadditive function. 
We say that $U$ is nonlattice if the group generated by its support
is dense in $\R$. We say that $z:\R_+\to\R_+$ is directly Riemann
integrable if $\int_{\R_+}\overline{z}_h(x)\mathrm dx<\infty$ for some $h>0$, and
$\int_{\R_+}(\overline{z}_h(x)-\underline{z}_h(x))\mathrm dx\to 0$ as $h\downarrow 0$,
where
$$\overline{z}_h(x)=\sup\left\{z(y):y\in \left[h\left\lfloor
      \frac{x}{h}\right\rfloor,h\left\lfloor\frac{x}{h}\right\rfloor+h\right]\right\}\,
,$$
and similarly for $\underline{z}_h$, with an $\inf$ instead of a
$\sup$. Note that these conditions imply that $z(t)\rightarrow 0$ as
$t \rightarrow \infty$. 
Let us recall the classical

\begin{lemma}[Key renewal theorem]
  \label{sec:potent-resolv-meas-1}
If $U$ is
nonlattice, then, for every directly Riemann integrable function 
$z:\R_+\to \R_+$, one has 
\begin{equation}
  \label{eq:31}
  z*U(t)=\int_{[0,t]}z(t-s)U(\mathrm ds)\underset{t\to\infty}{\longrightarrow}
  \frac{1}{\phi(0+)}\int_0^\infty z(s)\mathrm ds\, .
\end{equation} 
\end{lemma}

\noindent{\bf Proof. }
Observe that 
$$U+\delta_0=\sum_{n\geq 0}U_1^{*n}\, ,$$
so that $U+\delta_0$ is the renewal measure of the random walk with step
distribution $U_1$, so that this result is a consequence of
Blackwell's strong renewal theorem, see
\cite[Theorem V.4.3]{asmussen87}. There is a little subtlety here, since \cite{asmussen87} makes the working
assumption that the random walk step distribution does not charge
$\{0\}$. However, this is not a restriction, since, writing
$U_1=p\delta_0+(1-p)V_1$, with $V_1(\{0\})=0$, we have
$$\sum_{n\geq 0}U_1^{*n}=\sum_{n\geq
  0}\sum_{k=0}^n\binom{n}{k}((1-p)V_1)^{*k}p^{n-k}=\sum_{k\geq
  0}(1-p)^kV_1^{*k}\sum_{n\geq k}\binom{n}{k}p^{n-k}=\frac{\sum_{k\geq
  0}V_1^{*k}}{1-p}\, ,$$
so that 
$(1-p)(U+\delta_0)$ is the renewal measure of the random walk with step
distribution $V_1$, which does not charge $0$. \hfill$\square$

Some refinements of this result will be needed to obtain concentration
estimates in the subcritical and critical cases, but we postpone this discussion to
Section \ref{sec:subcritical-case}. 

\section{Analysis of the urn distribution function}
\label{sec:sizes-small-boxes}

Let $k\geq 2$ be a fixed integer. In order to apply Theorem \ref{thm:Karlin} to our situation, we need
to understand the behavior of $S_x:=S^{\mathbf{P}}_x$ as $x\downarrow
0$, where $\mathbf{P}$ is defined in \eqref{eq:3}, with $T=T^{(\alpha,\nu)}$.
By Proposition \ref{sec:part-repr-gene-1} and the discussion that
precedes it, 
the branchpoints
$b$ of
the tree $T^{(\alpha,\nu)}$ correspond exactly to the set $\mathcal{J}$ of times $t$ where an object of the
associated homogeneous fragmentation process $F^{(\nu)}$ splits into smaller
fragments. The total mass $\mu(T^{(\alpha,\nu)}_b)$ of the subtrees above $b$ equals the
size of the object before splitting, say $F^{(\nu)}_{i(t)}(t-)$ for some
$i(t)\geq 1$, and the measures
$\mu(T^{(\alpha,\nu)}_{b,i}),i\geq 1$ correspond to the sizes after splitting. These
can be written as $F^{(\nu)}_{i(t)}(t-)\Delta_j(t),j\geq 1$ for some sequence
$(\Delta_j(t),j\geq 1)\in \mathcal{S}$. In particular, we obtain
that 
\begin{equation}
  \label{eq:4}
  S_x=\sum_{t \in \mathcal J(F)} \sum_{i\geq 1} \mathbbm
  1_{\{i(t)=i\}}\mathbbm 1_{\left\{F^{(\nu)}_i(t-)^k \big(1-\sum_{j\geq
        1}\Delta^k_j(t)\big) \geq x \right\}}\, .
\end{equation}

In order to prove these results, we view $S_x=S_x(\infty)$ as the limiting value as
$t\to\infty$ of the adapted increasing process
\begin{equation}
  \label{eq:1}
  S_x(t)=\sum_{\substack{s \in \mathcal J(F)\\s\leq t}} \sum_{i\geq 1} \mathbbm 1_{\{i(s)=i\}}\mathbbm 1_{\left\{F^{(\nu)}_i(s-)^k \big(1-\sum_{j\geq 1}\Delta^k_j(s)\big) \geq x \right\}}.
\end{equation}
Since we are working with a homogeneous fragmentation, it holds \cite{berthfrag01} that the random measure
$$\sum_{t\in \mathcal{J}}\delta_{(t,i(t),(\Delta_j(t),j\geq 1))}$$
is a Poisson random measure on $\R_+\times \mathbb{N}\times
\mathcal{S}^\downarrow$ with intensity $\mathrm dt\mathbbm{1}_{\{t\geq 0\}}\#_{\mathbb{N}}(\mathrm di)\nu(\mathrm d\mathbf{s})$,
where $\#_{\mathbb{N}}$ is the counting measure on
$\mathbb{N}$. Therefore, the process $(S_x(t),t\geq 0)$ admits the compensator
\begin{equation}
  \label{eq:9}
    S^{(\mathrm p)}_x(t)
= \int_0^t \mathrm{d}s\sum_{i \geq 1}  
    f_k\left(\frac{x}{F^{(\nu)}_i(s)^k}\right)\, ,
 \end{equation}
  where we let $f_k:(0,\infty)\to \R$ be the nonincreasing function defined by 
  \begin{equation}
    \label{eq:10}
    f_k(x) = \nu\left(\sum_{i \geq 1} s_i^k \leq 1-x \right)\, ,\qquad
    x\in (0,1)
  \end{equation}
  and $f_k(x)=0$ for $x\geq 1$. This compensator is a nondecreasing
process, and we denote its limit as $t\to\infty$ by
$S^{(\mathrm{p})}_x=S^{(\mathrm{p})}_x(\infty)$.  
Note that the process $M_t=S_x(t)-S^{(\mathrm{p})}_x(t),t\geq 0$ is a local 
martingale, with quadratic variation 
$[M]_t=S_x(t),t\geq 0$. Since obviously $S_x(\infty)=S_x\leq 1/x$, we obtain that 
$M$ is in fact a true square-integrable martingale, and that 
$\mathbb{E}[S_x]=\mathbb{E}[S^{(\mathrm{p})}_x]$.  

The proofs of Theorem \ref{sec:thm-super} and Theorem \ref{sec:thm-sub-gammapos} will proceed in three main steps:
\begin{enumerate}[itemsep=0.1pt, topsep=0pt]
  \item Evaluate $\mathbb{E}[S^{(\mathrm{p})}_x]$ as $x\downarrow 0$. 
    \item Show that $S_x-S^{(\mathrm{p})}_x$ is small compared to 
    $\mathbb{E}[S^{(\mathrm{p})}_x]$. 
    \item Show that
      $S^{(\mathrm{p})}_x/\mathbb{E}[S^{(\mathrm{p})}_x]$ converges
      almost surely. 
\end{enumerate}

The first two points are easier and can be treated in an essentially
unified way in the supercritical or subcritical cases. The last point
is easy in the supercritical case, but much more delicate in the
subcritical case, so we treat these cases separately in Sections
\ref{sec:supercritical-case} and \ref{sec:subcritical-case}.

\subsection{Evaluation of $\mathbb{E}[S^{(\mathrm{p})}_x]$}

The first step consists in establishing the following, recalling the definitions of $C^{\mathrm{sub}}_\nu(k)$ and $C^{\mathrm{cr}}_\nu(k)$ in (\ref{eq:24}).

\begin{proposition}\label{sec:eval-mathbb}  Assuming \eqref{eq:hg},
  and, for $\gamma=0$, that $U$ is nonlattice, we have the following. 
\begin{enumerate}[itemsep=0pt, topsep=0pt]
\item[\emph{(i)}] If $k\gamma>1$, then
$$
\mathbb E\left[S^{(\mathrm p)}_x\right] ~\underset{x \downarrow 0} \sim ~ x^{-\gamma} \cdot \frac{c_\nu k^{\gamma}}{\phi(k\gamma-1)}.$$
\item[\emph{(ii)}] If $k\gamma<1$, then
$$
\mathbb E\left[S^{(\mathrm p)}_x\right] ~\underset{x \downarrow 0} \sim ~
 x^{-1/k} \cdot C^{\mathrm{sub}}_\nu(k).
$$
\item[\emph{(iii)}] If $k\gamma=1$, then 
$$
\mathbb E\left[S^{(\mathrm p)}_x\right] ~\underset{x \downarrow 0} \sim ~ 
  x^{-1/k} |\log(x)| \cdot C^{\mathrm{cr}}_\nu(k). 
$$
\end{enumerate}
\end{proposition}

To lighten a bit the notation, let us set
$g_k(x)=\mathbb{E}[S^{(\mathrm{p})}_x]$. By taking expectations in
\eqref{eq:9} with $t=\infty$, we obtain the formula
\begin{equation}
  \label{eq:5}
  g_k(x)=\int_{\mathbb R_+}
e^yf_k(xe^{ky}) U(\mathrm dy)=x^{-1/k}h_k*U\Big(\frac{1}{k}\log\Big(\frac{1}{x}\Big)\Big)\, ,
\end{equation}
where $U$ is as before the potential measure of the tagged fragment
subordinator $\xi_*$, and where
$h_k(y)=e^{-y}f_k(e^{-ky})\mathbbm{1}_{\{y\geq 0\}}$.
The function $g_k$ is
nonincreasing on $(0,1]$, with $g_k(1)=0$ and
\begin{equation}
  \label{eq:11}
  g_k(0+)=\nu(\mathcal{S})\int_{\mathbb R^+}
  e^yU(dy)=\infty\, ,
\end{equation}
by monotone convergence, because $U$ has infinite mass. 

We first record a simple result on the
asymptotic behavior of $f_k$.  

\begin{lemma}\label{sec:eval-mathbb-1}
Under \eqref{eq:hg}, it holds that
$$
f_k(x) \underset{x \downarrow 0}\sim  c_\nu \left( \frac{x}{k}\right)^{-\gamma }\quad \text{ and } \quad \sup_{x \in (0,1)} x^{\gamma} f_k(x) <\infty.
$$
\end{lemma}

\textbf{Proof.}
The second claim is an immediate consequence of the first one,
together with the fact that $f_k$ is nonincreasing on $(0,1)$. To
prove the asymptotic equivalent, we use that $\sum_{i\geq 2}s_i^k \leq
(1-s_1)^k$, so that 
$$
\nu\left(s_1^k +(1-s_1)^k \leq 1-x\right) \leq \nu\left(\sum_{i \geq
    1} s_i^k \leq 1-x \right) \leq \nu\left(s_1^k\leq 1-x\right)\, .
$$
By \eqref{eq:hg} and the fact that
$(1-x)^{1/k}=1-x/k+o(x)$ when $x \rightarrow 0$, we see that the upper
bound yields the right asymptotic equivalent. For the lower bound, we
observe that $s_1^k+(1-s_1)^k\leq 1-x$ if and only if $s_1\in
[a(x),1-a(x)]$, where $a(x)=x/k+o(x)$ as $x\downarrow
0$. Hence
$$\nu(s_1^k+(1-s_1)^k\leq 1-x)\geq \nu(s_1\leq 1-a(x))-\nu(s_1\leq
a(x))\, ,$$
and we conclude by \eqref{eq:hg} and the fact that
$\nu(s_1\leq a(x))\to 0$ as $x\to 0$, since $\nu(s_1\leq 1/2)<\infty$. 
$\hfill \square$

\vspace{0.1cm}

\noindent{\bf Proof of Proposition \ref{sec:eval-mathbb}. }
Statement (i) is an immediate application of the dominated convergence theorem,
 since, by the preceding lemma, $\sup_{x>0} x^{\gamma} e^{y}  f_k(xe^{ky})=O(e^{(1-k\gamma)y})$, and
 $\int_{\mathbb R^+} e^{(1-k\gamma)y}U(\mathrm dy)=1/\phi(k\gamma-1)<\infty$
 when $k\gamma>1$ by \eqref{eq:32}. 

For (ii), observe that $h_k(x)=O(e^{-(1-k\gamma)x})$, and that $h_k$
is continuous almost everywhere with respect to the Lebesgue measure. Hence, if
$k\gamma<1$, the function $h_k$ is directly Riemann integrable, and
the Key renewal theorem (Lemma \ref{sec:potent-resolv-meas-1}) implies
that, if $U$ is nonlattice, 
$$x^{1/k}g_k(x)=h_k*U\Big(\frac{1}{k}\log\Big(\frac{1}{x}\Big)\Big)\underset{x\downarrow
0}{\longrightarrow}\frac{1}{\phi(0+)}\int_0^\infty
h_k(x)\, \mathrm dx\, .$$

Let us finally prove (iii). Assuming $\gamma=1/k$, so in particular
$\gamma>0$, note that the L\'evy measure $\Xi$ of the subordinator $\xi_*$ is infinite and so $U$ is
necessarily nonlattice. We first write the
nonincreasing function 
$f_k$ as $f_k(x)=\zeta_k([x,1])$, where $\zeta_k$ is a nonnegative
measure. Then we have 
\begin{equation}
  \label{eq:21}
   g_k(x)=\int_{\mathbb R_+} e^y f_k(xe^{ky})U(\mathrm dy)= \int_{[x,1]} \zeta_k(\mathrm du) \int_{[0,\frac{1}{k}\log\left(\frac{u}{x}\right)]} e^y U(\mathrm dy).
\end{equation}
 By Lemma \ref{sec:potent-resolv-meas-1} applied to $z(t)=e^{-t}$, we
have that $\int_{[0,t]} e^y U(\mathrm dy)\sim
e^{t}/\phi'(0+)$ as $t\to\infty$. 
Thus, for $\delta>0$, there exists $a_{\delta} \in (0,\infty)$ such
that, for every $u \geq a_{\delta}x$, 
$$
(1-\delta) \frac{1}{\phi'(0+)} \frac{u^{1/k}}{x^{1/k}} \leq
\int_{[0,\frac{1}{k}\log\left(\frac{u}{x}\right)]} e^y U(\mathrm dy)
\leq (1+\delta) \frac{1}{\phi'(0+)} \frac{u^{1/k}}{x^{1/k}}.$$
This implies, since
 $f_k(x)\sim c_\nu k^{1/k}x^{-1/k}$ by Lemma \ref{sec:eval-mathbb-1}, 
 \begin{equation}
 \label{eq:HUupperbound}
  \int_{[x,1]} \zeta_k(\mathrm du)
  \int_{[0,\frac{1}{k}\log\left(\frac{u}{x}\right)]} e^y U(\mathrm dy)
  \leq O(x^{-1/k}) + (1+\delta) \frac{x^{-1/k}}{\phi'(0+)}
  \int_{[a_{\delta} x,1]} u^{1/k} \zeta_k(\mathrm du). 
  \end{equation}
  Integrating by parts and using again the asympotic behavior of
  $f_k$, we obtain
  $$\int_{[a_{\delta} x,1]} u^{1/k} \zeta_k(\mathrm du) =
  \frac{1}{k}\int_{a_\delta x}^1 u^{1/k-1}f_k(u)\, \mathrm du+O(1)=
  c_\nu k^{1/k-1}|\log(x)|+O(1),$$
  when $x\downarrow 0$. Together with (\ref{eq:HUupperbound}), this leads to 
$$
 \limsup_{x \rightarrow 0}  \frac{x^{1/k}}{|\log(x)|}g_k(x) \leq  (1+\delta) C^{\mathrm{cr}}_\nu(k)\, ,
 $$
for every $\delta>0$.  We argue similarly for the lower bound, concluding the proof of
 (iii). 
\hfill$\square$

{\remark When $\gamma=0$ and $U$ is lattice, then there can be
  oscillatory behavior for $g_k(x)$. However, it still holds that
  $g_k(x)=O(x^{-1/k})$ in this case, by an easy application of the
   renewal theorem in the lattice case. }

\subsection{Concentration of $S_x-S^{(\mathrm{p})}_x$}
\label{sec:concentration1}

We now turn to the property that $S_x$ and $S^{(\mathrm p)}_x$ are 
close. This comes from the following simple and very general variance estimate, 
which is valid even without assuming \eqref{eq:hg}. Compare with 
\cite[Lemma 2.9]{qshi15}, where it is shown that the inequality is in 
fact an equality, but we still give a proof for completeness. 

\begin{proposition}[Concentration]\label{sec:conc-s_xk-smathrmp_x}
For every $x>0$, it holds that 
$$\mathbb E\left[\left(S_x-S^{(\mathrm p)}_x \right)^2 \right]\leq g_k(x)$$
\end{proposition}

\noindent{\bf Proof.}
As discussed at the beginning of Section \ref{sec:sizes-small-boxes}, 
the process $M=(S_x(t)-S^{(\mathrm{p})}_x(t),t\geq 0)$ is a local 
martingale starting at $0$, with quadratic variation $(S_x(t),t\geq 0)$. Since 
$S_x(\cdot)$ is a pure-jump process with jumps of magnitude $1$ and 
$S^{(\mathrm{p})}_x(\cdot)$ is continuous, we can 
localize $M$ by a sequence of stopping times $T_n,n\geq 1$ with 
$T_n\to\infty$ such that $|M_{t\wedge T_n}|\leq n$, $S_x(t\wedge 
T_n)\leq n$ and $S_x^{(\mathrm{p})}(t\wedge T_n)\leq n$. Then we have by 
the stopping theorem 
$$\mathbb{E}[(M_{T_n})^2]=\mathbb{E}[(S_x(T_n))]=\mathbb{E}[(S_x^{(\mathrm{p})}(T_n))]$$
and Fatou's lemma implies the result. 
\hfill$\square$

\begin{corollary}
  \label{sec:conc-s_xk-smathrmp_x-1}
  It holds that 
  $$\frac{S_x-S^{(\mathrm{p})}_x}{g_k(x)}\overset{L^2}{\underset{x\downarrow  0}{\longrightarrow}}
0\, .$$
Moreover, assuming \eqref{eq:hg}, and, for $\gamma=0$, that $U$ is
nonlattice, and fixing $\lambda\in 
(0,1)$, the limit also holds almost surely 
along the values $x\in \{\lambda^n,n\geq 0\}$. 
\end{corollary}

\noindent{\bf Proof. }
We rewrite the statement of Proposition \ref{sec:conc-s_xk-smathrmp_x} in the form 
\begin{equation}
  \label{eq:6}
  \mathbb E\left[\left(\frac{S_x-S^{(\mathrm p)}_x}{g_k(x)} \right)^2 
  \right]\leq \frac{1}{g_k(x)}\, . 
\end{equation}
The statement on $L^2$ convergence then follows from this and \eqref{eq:11}. 
Now, assuming that \eqref{eq:hg} holds, Proposition 
\ref{sec:eval-mathbb} implies that the upper bound in \eqref{eq:6} is 
$O(x^{\max(\gamma,1/k)})$. In particular, it is summable over values of 
$x=\lambda^n,n\geq 0$, yielding the second claim. 
\hfill$\square$

Let us now assume for a minute that there exists a random variable 
$L$ such that, for every $\lambda\in (0,1)$, it 
holds that, almost surely, 
\begin{equation}
  \label{eq:7}
\frac{
  S^{(\mathrm{p})}_{\lambda^n}}{g_k(\lambda^n)}\underset{n\to\infty}{\longrightarrow}L\, . 
\end{equation}
Fixing $\lambda\in (0,1)$, and for a 
given $x\in (0,1]$, let $n$ be the unique integer 
      such that $\lambda^{n+1}<x\leq \lambda^n\, .$ By 
      monotonicity of $S_x$ and $g_k(x)$, 
      we have 
      $$\frac{S_{\lambda^{n}}-S^{(\mathrm{p})}_{\lambda^{n}}}{g_k(\lambda^{n+1})}
      +\frac{S^{(\mathrm{p})}_{\lambda^{n}}}{g_k(\lambda^{n+1})}\leq \frac{S_x}{g_k(x)}\leq 
      \frac{S_{\lambda^{n+1}}-S^{(\mathrm{p})}_{\lambda^{n+1}}}{g_k(\lambda^n)}
      +\frac{S^{(\mathrm{p})}_{\lambda^{n+1}}}{g_k(\lambda^n)}\, 
      .$$
      Assuming \eqref{eq:hg} and, when $\gamma=0$, that $U$ is nonlattice, and observing that, by 
      Proposition \ref{sec:eval-mathbb}, 
      $g_k(x)$ is regularly varying with 
      exponent $-\gamma'=-\max(\gamma,1/k)$, the second statement of Corollary 
      \ref{sec:conc-s_xk-smathrmp_x-1} implies that, almost surely, 
      $$ \lambda^{\gamma'}L\leq  \limsup_{x\downarrow 
        0}\frac{S_x}{g_k(x)}\leq \limsup_{x\downarrow 
        0}\frac{S_x}{g_k(x)}\leq \lambda^{-\gamma'}L\, .$$
      Since $\lambda\in (0,1)$ was arbitrary, we conclude that, almost surely 
      \begin{equation}
        \label{eq:8}
              \frac{S_x}{g_k(x)}\underset{x\downarrow0}{\longrightarrow}
              L\, . 
            \end{equation}
            
It remains to justify that \eqref{eq:7} holds, and identify the limit 
$L$, in order to apply Karlin's result, Theorem \ref{thm:Karlin}. The
supercritical case is easy, as we now discuss.

\subsection{The supercritical case}\label{sec:supercritical-case}

In this section, we prove Theorem \ref{sec:thm-super}.

\subsubsection{Almost sure convergence}

By Karlin's Theorem \ref{thm:Karlin}, the almost sure statements of
Theorem \ref{sec:thm-super} are direct corollaries of the following
proposition. 

\begin{proposition}
\label{prop:eq_S_x}
Assume \eqref{eq:hg} and $k\gamma>1$. Then, almost surely,
$$x^{\gamma} S_x \underset{x \downarrow 0} \rightarrow 
k^{\gamma}  \int_0^\infty c_\nu\sum_{i \geq 1}
F^{(\nu)}_i(t)^{k\gamma}\mathrm dt\, .$$
\end{proposition}

Note that the latter integral is indeed equal to $A^{(\nu)}_k$,
because of the time-change correspondence between the homogeneous fragmentation $\Pi^{(0,\nu)}$ and
the self-similar one $\Pi^{(\alpha,\nu)}$, as defined in Section \ref{sec:part-repr-gene}. Indeed, let us denote by
$\Pi^{(\alpha,\nu)}_i(t),i\geq 1$ the blocks of
$\Pi^{(\alpha,\nu)}(t)$ arranged in increasing order of their least
elements, and let $\tau_i^{\alpha}(t)=\tau^{(\alpha)}_{(j)}(t)$ for
any $j\in \Pi^{(\alpha,\nu)}_i(t)$, recalling \eqref{eq:28}. Then note
that $\Pi^{(\alpha,\nu)}_i(t)=\Pi^{(0,\nu)}_i(\tau^{(\alpha)}_{i}(t))$
for every $t\geq 0$ and $i\geq 1$, so that 
\begin{align*}
  \int_0^\infty \sum_{i \geq 1}
F^{(\nu)}_i(t)^{k\gamma}\mathrm dt&=\sum_{i\geq 1}\int_0^\infty
                                    |\Pi_i^{(0,\nu)}(t)|^{k\gamma}\mathrm dt\\
  &=\sum_{i\geq 1}\int_0^\infty
                                   |\Pi_i^{(0,\nu)}(t)| \mathrm
    d(\tau^{(1-k\gamma)}_i)^{-1}(t)\\
    &=\sum_{i\geq 1}\int_0^{\infty}
                                   |\Pi_i^{(1-k\gamma,\nu)}(t)| \mathrm
      dt=\int_0^{\infty}\sum_{i\geq 1}
                                   F_i^{(1-k\gamma,\nu)}(t)\mathrm
      dt\, .
\end{align*}

\noindent{\bf Proof of Proposition \ref{prop:eq_S_x}. }
Assume $\gamma>1/k$ and that \eqref{eq:hg} holds. One has, by 
(i) in Proposition \ref{sec:eval-mathbb}, 
$$\frac{S^{(\mathrm 
  p)}_x}{g_k(x)}\sim \frac{\phi(k\gamma-1)}{c_\nu k^\gamma}\int_0^\infty\sum_{i=1}^\infty x^\gamma 
  f_k\left(\frac{x}{F^{(\nu)}_i(t)^k}\right)\mathrm d t\, ,$$
  almost surely as $x\downarrow 0$. 
  For fixed $i,t$, one has by Lemma \ref{sec:eval-mathbb-1}
  \begin{equation}
    \label{eq:36}
      x^\gamma 
  f_k\left(\frac{x}{F^{(\nu)}_i(t)^k}\right)\longrightarrow c_\nu k^\gamma 
  F^{(\nu)}_i(t)^{k\gamma}\, ,
\end{equation}
while being dominated by $F^{(\nu)}_i(t)^{k\gamma}\sup_{y\in (0,1)}y^\gamma
  f_k(y)$. 
  Since $k\gamma>1$, we have 
  $$\mathbb{E}\left[\int_0^\infty\sum_{i\geq 1}F^{(\nu)}_i(t)^{k\gamma} \mathrm dt\right]=\int_0^\infty\mathbb{E}\big[F^{(\nu)}_*(t)^{k\gamma-1}\big]\mathrm d t=\frac{1}{\phi(k\gamma-1)}<\infty\, 
,$$
where $F_*$ is the tagged fragment. 
Therefore, the dominated convergence theorem applies and gives, almost 
surely, 
\begin{equation}
  \label{eq:37}
  \frac{S^{(\mathrm{p})}_{x}}{g_k(x)}\underset{x\downarrow 
0}{\longrightarrow}\phi(k\gamma-1)\int_0^\infty\sum_{i\geq 
1}F^{(\nu)}_i(t)^{k\gamma}\, \mathrm d t\, ,
\end{equation}
and this obviously implies \eqref{eq:7} and identifies $L$. We 
conclude that \eqref{eq:8} holds, and Proposition \ref{prop:eq_S_x} follows by using again 
the asymptotic behavior of $g_k(x)$ from 
(i) in Proposition \ref{sec:eval-mathbb}. 
\hfill$\square$

\subsubsection{$L^2$ convergence}\label{sec:l2-convergence}

We still assume \eqref{eq:hg} and $k\gamma>1$, and now prove the statement on $L^2$ convergence in Theorem
\ref{sec:thm-super}.
To that end, recalling that $\mathbf{P}$ denotes the urn sizes, we note that
$\mathbb E[N^{(\nu)}_n(k)\, |\, \mathbf{P}]=n\int_0^1(1-x)^{n-1}S_x\mathrm dx$
(see \cite[Equation (4)]{GPY06}). Denoting this quantity by
$\widetilde{N}^{(\nu)}_n(k)$, we have, by \cite[Equation (60)]{karlin67}, 
$$\mathbb E[(N^{(\nu)}_n(k)-\widetilde{N}^{(\nu)}_n(k))^2\, |\, \mathbf{P}]\leq
\widetilde{N}^{(\nu)}_n(k)\, ,$$
so that it suffices to show that $n^{-\gamma}\widetilde{N}^{(\nu)}_n(k)$
converges in $L^2$ to the a.s.\ limit $\Gamma(1-\gamma)L=\lim n^{-\gamma}N^{(\nu)}_n(k)$
obtained in the previous paragraph, to conclude that
$n^{-\gamma}N^{(\nu)}_n(k)$ also converges to $\Gamma(1-\gamma)L$ in $L^2$. 

Denoting by $B(a,b)$ the Beta function, we obtain, after some elementary manipulations
$$
  \frac{B(1,n)}{B(1-\gamma,n)}\widetilde{N}^{(\nu)}_n(k)-L=
  \int_0^1 (x^\gamma S_x-L)\frac{x^{-\gamma}(1-x)^{n-1}\mathrm
    dx}{B(1-\gamma,n)}\, ,$$
  and so, by Jensen's inequality,
  \begin{equation*}
    \label{eq:35}
   \mathbb E\left[
     \left(\frac{B(1,n)}{B(1-\gamma,n)}\widetilde{N}^{(\nu)}_n(k)-L\right)^2\right]\leq
   \frac{n^{\gamma-1}}{B(1-\gamma,n)}\int_0^n
   \mathbb{E}\left[\left(\left(\frac{x}{n}\right)^\gamma S_{\frac{x}{n}}-L\right)^2\right]
 x^{-\gamma}\left(1-\frac{x}{n}\right)^{n-1}\mathrm
    dx.
  \end{equation*}
We see that in turn that the wanted convergence will be a consequence of
the fact that $\mathbb E[(x^\gamma S_x-L)^2]$ converges to $0$ as
$x\downarrow 0$, by an immediate application of dominated
convergence. From Corollary \ref{sec:conc-s_xk-smathrmp_x-1} and
Proposition \ref{sec:eval-mathbb} (i), it suffices to show this
convergence for $S^{(\mathrm{p})}_x$ in place of $S_x$. 

By \eqref{eq:37}, we already know that $x^\gamma S_x^{(\mathrm{p})}$
converges a.s.\ to $L$, so it suffices to show that it is also bounded
in $L^q$ for some $q>2$. However, we have, by the remark just after
\eqref{eq:36}, 
$$
  x^{\gamma}S_x^{(\mathrm{p})}\leq C\left(\int_0^1\mathrm dt \sum_{i\geq
  1}F^{(\nu)}_i(t)^{k\gamma}+\int_1^\infty\frac{\mathrm dt}{t^2}\sum_{i\geq
  1}F^{(\nu)}_i(t) t^2F^{(\nu)}_i(t)^{k\gamma-1}\right)\, ,
$$
for some $C\in (0,\infty)$. 
Now, using the fact that $\sum_{i\geq 1}F_i(t)=1$, and Jensen's
inequality (note that $t^{-2}\mathrm dt\sum_{i\geq 1}F^{(\nu)}_i(t)\delta_i$
is a probability measure), we obtain
\begin{align*}
  \mathbb{E}[(x^\gamma S_x^{(\mathrm{p})})^q]
&\leq 2^{q-1}C\left(1+\int_1^\infty \mathrm dt\,
  t^{2q-2}\mathbb{E}\left[\sum_{i\geq
                                                1}F^{(\nu)}_i(t)^{q(k\gamma-1)+1}\right]\right)\\
  &= 2^{q-1}C\left(1+\int_1^\infty\mathrm dt\,
    t^{2q-2}e^{-t\phi(q(k\gamma-1))}\right)\, ,
\end{align*}
using the many-to-one formula in the last step and the definition of
the Laplace exponent $\phi$ of the tagged fragment subordinator. The
latter quantity is finite for every $q>0$, yielding the result.

\section{The subcritical and critical cases}\label{sec:subcritical-case}

The goal of this section is to prove the following proposition and then Theorem \ref{sec:thm-sub-gammapos}.

\begin{proposition}
  \label{sec:sizes-small-boxes-1}
  Assume  \eqref{eq:hg}, \eqref{eq:exp}, and also
  \eqref{eq:gamma0} if $\gamma=0$. Then, almost surely,
  \begin{enumerate}[itemsep=0pt, topsep=0pt]
\item[\emph{(i)}] if $k\gamma=1$, 
$$
S_x ~\underset{x \downarrow 0} \sim~ C^{\mathrm cr}_{\nu}(k)\cdot x^{-\frac{1}{k}} |\log(x)|,
$$
\item[\emph{(ii)}] if $k\gamma<1$, 
$$
S_x ~\underset{x \downarrow 0} \sim~ C^{\mathrm{sub}}_\nu(k)\cdot x^{-1/k}\, .
$$
\end{enumerate}
\end{proposition}

We will prove this using \eqref{eq:7}.  Since the limits are now
deterministic, this requires the following variance estimate.

\begin{proposition}
  \label{sec:conc-comp-1}
  Assuming \eqref{eq:hg}, \eqref{eq:exp} and $k\gamma\leq 1$, and also
  \eqref{eq:gamma0} if $\gamma=0$, it holds that for some $C\in(0,\infty)$, 
  $$\mathrm{Var}(S^{(\mathrm{p})}_x)\leq
  C\frac{g_k(x)^2}{|\log(x)|^2}\, .$$
\end{proposition}

In fact, for $k\gamma<1$, the logarithm in the denominator can be improved to a negative
power of $x$. Given this statement, we can finish the proof of
Proposition \ref{sec:sizes-small-boxes-1}.

\noindent{\bf Proof of Proposition \ref{sec:sizes-small-boxes-1}. }
Assuming
\eqref{eq:hg} and \eqref{eq:exp}, Proposition \ref{sec:conc-comp-1}
implies, for every $\lambda\in (0,1)$,  
    $$\mathrm{Var}\left(\frac{S^{(\mathrm{p})}_{\lambda^n}}{g_k(\lambda^n)}\right)\leq 
      \frac{C}{n^2\log(\lambda)^2}\ , $$
      and this is summable in $n$. Consequently, almost surely, 
      $$\frac{S^{(\mathrm{p})}_{\lambda^n}}{g_k(\lambda^n)}\underset{n\to\infty}{\longrightarrow }
      1\, .$$ Hence, \eqref{eq:7}, and, therefore, \eqref{eq:8}, hold
      with $L=1$. The a.s.\ convergence statement of Proposition \ref{sec:sizes-small-boxes-1} follows by
      using the asymptotic behavior of
      $g_k(x)$ given by (ii) and (iii) in
      Proposition \ref{sec:eval-mathbb}.
\hfill$\square$

This immediately leads to Theorem \ref{sec:thm-sub-gammapos}:

 \noindent{\bf Proof of Theorem \ref{sec:thm-sub-gammapos}. }     
Proposition \ref{sec:sizes-small-boxes-1} directly implies the almost
sure convergence statements of Theorem 
\ref{sec:thm-sub-gammapos} by a use of Karlin's result. 

To obtain the statement about convergence in $L^2$, we can repeat {\em verbatim} the
first part of the argument of Section \ref{sec:l2-convergence}, replacing the exponent
$\gamma$ by $1/k$ everywhere, and adding a $\log$ factor in the case where $k
\gamma=1$. 
We obtain that the $L^2$ convergence is a consequence of
the fact that $x^{1/k} S_x^{(\mathrm{p})}$ (resp.\ $|\log(x)|^{-1}x^{1/k} S_x^{(\mathrm{p})}$) converges in $L^2$ to its
almost sure limit when $k\gamma<1$ (resp.\ when $k\gamma=1$). But this is immediate by Proposition
\ref{sec:conc-comp-1} and Proposition \ref{sec:eval-mathbb} (ii) and (iii). \hfill$\square$

It remains to prove Proposition \ref{sec:conc-comp-1}. This will be
done in Section \ref{sec:conc-comp}, after we gather some refined results in the next section on the potential and resolvent measures
$U$ and $U_1$ of the tagged fragment subordinator $\xi_*$. 

\subsection{Refined renewal estimates}\label{sec:refin-renew-estim}

The following lemmas provide some regularity results for the measures
$U_1$ and $U$, under our working assumptions. 

\begin{lemma}
  \label{sec:fragmentation-trees-1}
  \emph{(i)} If \eqref{eq:hg} holds 
  with $\gamma\in (0,1)$, then, for every $t> 0$, the law of 
  $\xi_*(t)$ admits a density that is infinitely differentiable and 
  has bounded derivatives of all orders. In particular, $U_1$ is
  absolutely continuous. 

\emph{(ii)}  Assume that \eqref{eq:hg} holds 
  with $\gamma=0$, and that \eqref{eq:gamma0} holds. 
  Then for every $t>0$, the singular part of the law of $\xi_*(t)$ is
  $e^{-t\Xi(\R_+)}\delta_0$. In particular, $U_1$ is absolutely
  continuous on $(0,\infty)$.   
  
  \emph{(iii)} If \eqref{eq:exp} holds, then there exists $\varepsilon>0$ such
  that $\int_0^\infty U_1(\mathrm  dx)e^{\varepsilon x}<\infty$.  
\end{lemma}

\noindent{\bf Proof. }
Assume \eqref{eq:hg} with $\gamma>0$, so that $2-\gamma\in (0,2)$. Then, by a criterion of Orey, see \cite[Theorem 28.3]{sato13}, (i) is a 
consequence of the fact that 
\begin{equation}
  \label{eq:22}
  \liminf_{r\downarrow 0}\frac{1}{r^{2-\gamma}}\int_0^rx^2\, \Xi(\mathrm dx)>0\, . 
\end{equation}
Let us prove this fact. We write, for $r<\log(2)$, 
\begin{align*}
  \int_0^rx^2\, \Xi(\mathrm dx)&=\sum_{i\geq 
                                 1}\int_{e^{-r}}^1(-\log(y))^2\, y\, 
                                 \nu(s_i\in \mathrm dy)\\
                                 &=\int_{e^{-r}}^1(-\log(y))^2\, y\, 
                                   \nu(s_1\in \mathrm dy)\, , 
\end{align*}
since we have $s_i \leq 1/2$ for every $i\geq 2$, $\nu(\mathrm 
d\mathbf{s})$-almost surely. In turn, this last expression is 
$$
   \geq 2\log(2)^2\int_{e^{-r}}^1(1-y)^2\nu(s_1\in \mathrm dy)\geq \log(2)^2\int_0^{1-e^{-r}}y\mathrm dy\big(\nu(s_1\leq 
   1-y)-\nu(s_1\leq e^{-r})\big)\, ,$$
   and by \eqref{eq:hg}, we see that $\nu(s_1\leq 
   e^{-r})\sim c_\nu r^{-\gamma}$ and $\nu(s_1\leq 1-y)\sim c_\nu 
   y^{-\gamma}$, so that \eqref{eq:22} holds. 

Now assume that \eqref{eq:hg} holds with $\gamma=0$, and that
$\sum_{i\geq 1}\nu(s_i\in \mathrm dx)$ is absolutely continuous. As
already observed, the same is true of $\Xi(\mathrm dx)$. Since the latter is
finite, $\xi_*$ is a compound Poisson process, and the law of
$\xi_*(t)$ is
$$e^{-t\Xi(\R_+)}\delta_0+e^{-t\Xi(\R_+)}\sum_{m\geq 1}\frac{t^m}{m!}\Xi^{*m}\, ,$$
where the last sum is absolutely continuous. This proves
the first case of (ii). The second case is obvious.

   To prove (iii), we use the fact that the Laplace transform of $U_1$ is 
   given by $\mathcal{L}_1(q)=(1+\phi(q))^{-1}$. 
If \eqref{eq:exp} holds, then, as observed around \eqref{eq:27}, $\phi$ can be analytically continuated
on some interval $(-\eta,\infty)$, so that 
   $\mathcal{L}_1$ can also be continued in $(-\varepsilon,\infty )$, where 
   $\varepsilon=\eta\wedge\inf \{q>0:\phi(-q)=-1\}$. Since the 
   Taylor coefficients of $\mathcal{L}_1$ are given by the moments of 
   $U_1$ (with alternating signs), we conclude that the expression 
   $\mathcal{L}_1(q)=\int_0^\infty e^{-q x}U_1(\mathrm dx)$
   remains valid in this whole domain. 
\hfill$\square$

In the following statement, we let $Z(t)=\int_{[0,t]} e^{-(t-y)}U(\mathrm 
dy)$ be the convolution of the function $z:x\mapsto
e^{-x}\mathbbm{1}_{\{x\geq 0\}}$ with the 
potential measure $U$. 

\begin{lemma}\label{sec:refin-renew-estim-1}
  Assume \eqref{eq:hg} and \eqref{eq:exp} hold with $\gamma\in
  [0,1)$. If $\gamma=0$, we further assume \eqref{eq:gamma0}. 
  Then
  there exists $C\in (0,\infty)$ and $\epsilon\in (0,1)$ such that, 
  for every $t,s\geq 0$ with $s\leq t$, 
  \begin{equation}
    \label{eq:26}
    |Z(t)-Z(t-s)|\leq C(U([0,s])\wedge 1)e^{\epsilon(s-t)}\, .
  \end{equation}
\end{lemma}

\noindent{\bf Proof. }
By Lemma
\ref{sec:fragmentation-trees-1} (i) and (ii), it holds that 
$U_1$ can be written as $p\delta_0+(1-p)V_1$, where $p\in [0,1)$ and
$V_1(\mathrm dx)=v_1(x)\mathrm dx$ is absolutely continuous. Moreover,
as observed in the proof of Lemma \ref{sec:potent-resolv-meas-1}, 
$(1-p)(U+\delta_0)$ is the renewal measure of the random walk with step
distribution $V_1$. In particular, $U$ is absolutely continuous on $(0,\infty)$,  with
density $u=(1-p)^{-1}\sum_{n\geq 1}v_1^{*n}$. 

Moreover,
$V_1$ admits small
exponential moments by Lemma \ref{sec:fragmentation-trees-1} (iii), and has (necessarily finite)
mean $\phi'(0+)/(1-p)$ by \eqref{eq:25}. 

In particular, \cite[Corollary
VII.1.3]{asmussen87} implies that $Z$ is a bounded function. Moreover, 
by (ii) and (iii) in 
\cite[Theorem 
VII.2.10]{asmussen87}, there exists $\epsilon\in (0,1)$ such that,
as 
$t\to\infty$, 
\begin{equation}
  \label{eq:30}
  u(t)=\frac{1}{\phi'(0+)}+O(e^{-\epsilon t}) \qquad \mbox{ and }\qquad
  Z(t)=\frac{1}{\phi'(0+)}+O(e^{-\epsilon t})\, .
\end{equation}
Now, we observe that the function $Z$ is differentiable, with
$$Z'(t)=-Z(t)+u(t)=O(e^{-\epsilon t})\, ,$$
as $t\to\infty$, and therefore, there exists $c\in (0,\infty)$ such
that $|Z'(t)|\leq ce^{-\epsilon t}$ for every
$t\geq 1$. By integrating this bound, we obtain that for $t,s\geq 0$ with $t-s\geq 1$, 
$$|Z(t)-Z(t-s)|\leq \frac{c}{\epsilon}(1-e^{-\epsilon
s}) e^{-\epsilon (t-s)}\, ,$$
which is an inequality of the wanted form, since $(1-e^{-\epsilon
s})\leq s\wedge 1\leq 
C(U([0,s])\wedge 1)$ for some finite constant $C$, by subadditivity of
$a\mapsto U([0,a])$. 
It remains to discuss the situation where $0\leq s\leq t$ with
$t-s\leq 1$. 
Let us first assume that $0< s\leq t\leq 2$, and observe that
$$Z(t)-Z(t-s)\leq \int_{[t-s,t]}U(\mathrm dy)\leq U([0,s])\, ,$$
by
subadditivity. Moreover, we have 
$$Z(t)-Z(t-s)\geq
Z(t)-(1+4s)Z(t)\geq
-4s U([0,2])\geq -8U([0,s])\, ,$$
so that \eqref{eq:26} holds in this case. Finally, for $t\geq 2$ and
$t-s\leq 1$, necessarily $s\geq 1$ and the result is an immediate consequence of the boundedness
of $Z$.
\hfill $\square$

\subsection{Concentration of $S^{(\mathrm{p})}_x$: proof of
  Proposition \ref{sec:conc-comp-1}}\label{sec:conc-comp}

We use some martingale
concentration techniques. Fix $x>0$, and consider the martingale
$M'_t=\mathbb E[S^{(\mathrm p)}_x\, |\, \mathcal F_t],t\geq 0$, where $(\mathcal F_t,t\geq 0)$ denotes the natural filtration of the process $F^{(\nu)}$. Note that
$M'_0=g_k(x)$, while
$M'_\infty=S^{(\mathrm p)}_x$.
Therefore, one has
$$\mathbb E[|S_x^{(\mathrm{p})}-g_k(x)|^2]\leq \mathbb E[[M',M']_\infty]$$
and our task is to control the
quadratic variation of the martingale $M'_t$. 
By
using the fragmentation property \cite{berthfrag01} representing the process
$(F^{(\nu)}(t+s),s\geq 0)$ as the superimposition of processes
$(F^{(\nu)}_i(t)F^{(\nu,i)}(s),s\geq 0),i\geq 1$, where the processes $F^{(\nu,i)}$
are i.i.d.\ copies of $F^{(\nu)}$, we have 
\begin{align*}
  M'_t&=\int_0^t \mathrm ds \sum_{i\geq 1}f_k\left(\frac{x}{F^{(\nu)}_i(s)^k}\right)+\sum_{i\geq 1}\mathbb E\left[\sum_{j\geq 1}\int_0^\infty f_k\left(\frac{x}{F^{(\nu)}_i(t)^kF^{(\nu,i)}_j(s)^k}\right)\, \Big|\,
  F_i(t)\right]\\
  &=\int_0^t \mathrm ds \sum_{i\geq 1}f_k\left(\frac{x}{F^{(\nu)}_i(s)^k}\right)+\sum_{i\geq
    1}g_k\left(\frac{x}{F^{(\nu)}_i(t)^k}\right)\, .
\end{align*}
Note that the sum is really a finite sum for every $t\geq 0$, since
$g_k(x)=0$ for $x\geq 1$. 
From this, it is easy to see that the martingale $M'-M'_0$ is of
finite variation, and hence is purely discontinuous, with quadratic
variation equal to the sum of the squares of its jumps, that is, 
$$[M',M']_\infty=\sum_{t\in\mathcal{J}(F)}\sum_{i\geq 1}\mathbbm
1_{\{i(t)=i\}}\left(\left(\sum_{j\geq
      1}g_k\left(\frac{x}{F^{(\nu)}_i(t-)^k\Delta_j(t)^k}\right)\right)-g_k\left(\frac{x}{F^{(\nu)}_i(t-)^k}\right)\right)^2\,
.$$
Taking expectations and using again a compensation formula, and then the
many-to-one formula, we obtain
\begin{align*}
  \mathbb E\left[[M',M']_\infty\right]&=\mathbb E\left[\int_0^\infty
  \mathrm dt \sum_{i\geq 1}\int_{\mathcal S} \nu(\mathrm d\mathbf s) \left(\left(\sum_{j\geq
  1}g_k\left(\frac{x}{F^{(\nu)}_i(t)^ks_j^k}\right)\right)-g_k\left(\frac{x}{F^{(\nu)}_i(t)^k}\right)\right)^2\right]\\
  &=\int_0^\infty U(\mathrm d y)e^y \int_{\mathcal S} \nu(\mathrm d\mathbf s) \left(\left(\sum_{j\geq
  1}g_k\left(\frac{xe^{ky}}{s_j^k}\right)\right)-g_k\left(xe^{ky}\right)\right)^2\, .
\end{align*}
At this point, we need the following technical estimate, whose proof
is postponed to after the current discussion. 

\begin{lemma}\label{sec:proof-prop-refs}
Assume \eqref{eq:hg}, \eqref{eq:exp}, and also
  \eqref{eq:gamma0} if $\gamma=0$, and let $k \in \mathbb N$. 
\vspace{-\topsep}  
\begin{itemize}
\setlength{\itemsep}{0.1pt}
\item If $k\gamma <1$, then there exists $\varepsilon>0$ such that
  \begin{equation}
    \label{eq:17}
    \int_{\mathcal S} \nu(\mathrm d \mathbf s) \left(\sum_{i\geq 1} g_k(xs_i^{-k})-g_k(x) \right)^2=O \left(x^{-\frac{2}{k}+2\varepsilon} \right).
  \end{equation}
\item  If $k\gamma =1$, then
  \begin{equation}
    \label{eq:18}
    \int_{\mathcal S} \nu(\mathrm d \mathbf s) \left(\sum_{i\geq 1} g_k(xs_i^{-k})-g_k(x) \right)^2=O \left(x^{-\frac{2}{k}} \right).
  \end{equation}
\end{itemize}
\end{lemma}

If $k\gamma<1$, Lemma \ref{sec:proof-prop-refs} yields
$$ \mathbb E\left[[M',M']_\infty\right]=O(
x^{-2/k+2\varepsilon}) \int_0^\infty U(\mathrm d y)e^{-y(1-2k\varepsilon)} \, ,$$
so that $\mathbb E\left[[M',M']_\infty\right]=
O(x^{2\varepsilon}g_k(x)^2)$, 
if $\varepsilon>0$ is chosen small enough so that
$1-2k\varepsilon>0$, since the last displayed integral then converges, and
since $x^{-1/k}=O(g_k(x))$ 
by (ii) in Proposition \ref{sec:eval-mathbb}. If $k\gamma=1$, on the other hand, Lemma  \ref{sec:proof-prop-refs} gives
$$ \mathbb E\left[[M',M']_\infty\right]=O(
x^{-2/k})\int_0^\infty U(\mathrm d y)e^{-y} \, .$$
Since $x^{-1/k}=O(g_k(x)/|\log(x)|)$ by (iii) in Proposition
\ref{sec:eval-mathbb}, we conclude that $\mathbb E\left[[M',M']_\infty\right]\leq
O(g_k(x)^2/|\log(x)|^2)$. The latter bound is thus valid for every $k,\gamma$
such that $k\gamma\leq 1$, as wanted. This concludes the proof of
Proposition \ref{sec:conc-comp-1}, except for
Lemma \ref{sec:proof-prop-refs}. 

\subsection{Proof of Lemma \ref{sec:proof-prop-refs} }

We start with a technical lemma. 

\begin{lemma}
\label{lem:increments_g}
Assume \eqref{eq:hg} and \eqref{eq:exp}, as well as \eqref{eq:gamma0}
if $\gamma=0$. Then for all $k\in \mathbb N$, there exists $\varepsilon>0$ such that:
\vspace{-\topsep}
\begin{itemize}
\setlength{\itemsep}{0.1pt}
\item If $\gamma=0$, there exists $\kappa \in (0,\infty)$ such that
  for all $\lambda \in (0,1)$ 
  and all $x \in (0,1)$
  \begin{equation}
    \label{eq:14}
      \left|g_k(x)-\lambda^{-1} g_k(\lambda^{-k} x) \right| \leq
      \kappa x^{-\frac{1}{k}+\varepsilon} \cdot \lambda^{-k\varepsilon
      }.
      \end{equation}
  \item If $\gamma \in (0,1)$ and $k\gamma <1$,  there exists $\kappa
    \in (0,\infty)$ such that for all $\lambda \in (0,1)$
    and all $x \in (0,1)$
    \begin{equation}
      \label{eq:15}
          \left|g_k(x)-\lambda^{-1} g_k(\lambda^{-k} x) \right| \leq
          \kappa x^{-\frac{1}{k}+\varepsilon}
          \left(\lambda^{-k\varepsilon } \mathbbm 1_{\{\lambda \leq
              1/2\}}+(1-\lambda)^{\gamma} \mathbbm 1_{\{\lambda \geq  1/2\}}
          \right).
        \end{equation}
      \item If $k\gamma =1$, for all $\lambda \in (0,1)$
        \begin{equation}
          \label{eq:16}
                  \left|g_k(x)-\lambda^{-1} g_k(\lambda^{-k} x)
                  \right| \leq \kappa x^{-\frac{1}{k}}
                  \left(|\log(\lambda)| \mathbbm 1_{\{\lambda \leq
                      1/2\}}+(1-\lambda)^{\gamma}\mathbbm 1_{\{\lambda\geq
                      1/2\}} \right).
                \end{equation}
              \end{itemize}
\end{lemma}

\textbf{Proof.}  Let $\lambda,x \in (0,1)$. 

{\bf We first assume that $\lambda^{-k} x \geq 1$.} In this case, we note that  
$$\left|g_k(x)-\lambda^{-1} g_k(\lambda^{-k} x) \right| =g_k(x),$$
which is either $O(x^{-1/k})$ when
$k\gamma<1$, or $O(x^{-1/k}|\log(x)|)$ when $k\gamma=1$, by
Proposition \ref{sec:eval-mathbb}. This yields the wanted bound
\eqref{eq:14}, and also, if we further assume that $\lambda\leq 1/2$,
the bounds \eqref{eq:15} and \eqref{eq:16}. 
When $\lambda>1/2$, observe that $x>1/2^k$ by our initial
assumption, so that 
\begin{eqnarray*}
g_k(x)&=&\int_{[0,\log(x^{-1/k})]} e^y f_k(xe^{ky})U(\mathrm dy) \\
&\leq& C U([0,\log(x^{-1/k})])
\end{eqnarray*}
for some finite constant $C>0$.
By \eqref{eq:13} and Proposition 1.5 of \cite{BertoinSubSF}, we have 
\begin{equation}
  \label{eq:12}
  \Gamma(1+\gamma)U([0,a])\underset{a\downarrow 0}{\sim} \frac{a^{\gamma}}{c_\nu\Gamma(1-\gamma)}\, ,
\end{equation}
and therefore,
$U([0,\log(x^{-1/k})])=O((1-x)^\gamma)$ as $x\to 1$. 
Recalling our initial assumption that $\lambda^{-k} x \geq
1$, we obtain that there exists a finite constant  $C''>0$ such that, 
when $\lambda > 1/2$, 
$$g_k(x)\leq C''(1-\lambda)^{\gamma}\, .$$
Since $x>1/2^k$ under our working
assumptions, this yields \eqref{eq:15} and 
\eqref{eq:16}. Observe that, so far, we can choose the value of
$\varepsilon$ arbitrarily. 

{\bf Now we assume that $\lambda^{-k} x <1$.} Recalling
\eqref{eq:21}, we have that 
\begin{eqnarray}
\notag && \left|g_k(x)-\lambda^{-1} g_k(\lambda^{-k} x) \right| \\ \label{Int1}
&\leq& \int_{[x,\lambda^{-k} x]} \zeta_k(\mathrm du)  \int_{[0,\frac{1}{k}\log(\frac{u}{x})]} e^yU(\mathrm dy) \\ \label{Int2}
       &+& \int_{[\lambda^{-k} x,1]} \zeta_k(\mathrm du) \left(\frac{u}{x}\right)^{1/k}\left|
Z\left(\frac{1}{k}\log\left(\frac{u}{x}\right)\right)
        -
           Z\left(\frac{1}{k}\log\left(\frac{u}{x}\right)+\log(\lambda)\right)  \right|\, ,
\end{eqnarray}
where $Z$ is defined before Lemma \ref{sec:refin-renew-estim-1}. 
We will bound separately the two integrals
(\ref{Int1}) and (\ref{Int2}), starting with the former. By using again \eqref{eq:12}, we obtain that, for $a\in [0,1]$, 
\begin{equation*}
  \label{eq:19}
  \int_{[0,a]}e^yU(\mathrm dy) =O(a^{\gamma})\, .
\end{equation*}
Now note that, in the integrand of
\eqref{Int1}, we have  $\log(u/x)/k\leq |\log(\lambda)|$. Hence, when
$\lambda> 1/2$, we have $\log(u/x)/k\leq |\log(2)|\leq 1$, and so we can
apply the above estimate and obtain that (\ref{Int1}) is bounded by a
multiple of $|\log(\lambda)|^{\gamma}\int_{[x,\lambda^{-k} x]}
\zeta_k(\mathrm du)$, which itself is smaller than a multiple of
$x^{-\gamma}(1-\lambda)^{\gamma}$, by Lemma
\ref{sec:eval-mathbb-1}. This is sufficient for the purpose of
(\ref{eq:14}--\ref{eq:16}). 

When $\lambda \leq 1/2$, we use the right-hand side of \eqref{eq:30} to
bound the integral (\ref{Int1}) by a multiple of 
$\int_{[x,\lambda^{-k} x]} (u/x)^{1/k} \zeta_k(\mathrm du).$ Integrating by parts, we see that 
\begin{align*}
    \int_{[x,\lambda^{-k} x]} (u/x)^{1/k} \zeta_k(\mathrm du)&=
 \zeta_k([x,x\lambda^{-k}]) 
+ \frac{1}{k}\int_1^{\lambda^{-k}} v^{1/k-1}
\zeta_k([xv,x\lambda^{-k}]) \mathrm dv\\
  &\leq f_k(x)+\frac{1}{k}\int_1^{\lambda^{-k}}v^{1/k-1}f_k(xv)\mathrm
    dv\, .
\end{align*}
By \eqref{eq:hg}, this is bounded above by a multiple of
$$x^{-\gamma}\left(1+\int_1^{\lambda^{-k}}v^{\frac{1}{k}-\gamma-1}\mathrm
  dv\right)=\left\{\begin{array}{lcl}
                     O(\lambda^{-1+k\gamma}x^{-\gamma}) & \text{ if } &k\gamma<1\\
                     O(|\log(\lambda)|x^{-1/k}) & \text{ if }&
                                                              k\gamma=1\, .
                   \end{array}
                 \right. $$
This is enough to get a bound compatible with (\ref{eq:14}--\ref{eq:16}): indeed, it is immediate when $k\gamma=1$, while when $k\gamma<1$, 
provided that $\varepsilon<1/k-\gamma$, we have that 
$$\lambda^{-1+k\gamma}=\lambda^{-k\varepsilon-k(\frac{1}{k}-\gamma-\varepsilon)} \leq \lambda^{-k\varepsilon} x^{-\frac{1}{k}+\varepsilon+\gamma}$$ since $\lambda^{-k} x <1$. 

It remains to bound from above the integral (\ref{Int2}).
By Lemma \ref{sec:refin-renew-estim-1} and \eqref{eq:12}, it is smaller than a multiple of 
$$
 \min(|\log(\lambda)|,1)^{\gamma} \cdot  \lambda ^{-k\varepsilon} \cdot x^{-\frac{1-k\varepsilon}{k}} \cdot \int_{[\lambda^{-k} x,1]} \zeta_k(\mathrm du) u^{\frac{1-k\varepsilon}{k}}.
$$
for some $\varepsilon>0$ small enough.
When $k\gamma<1$, then $\int_{(0,1]}
\zeta_k(\mathrm du) u^{\frac{1-k\varepsilon}{k}}<\infty$ by \eqref{eq:hg}, and this yields \eqref{eq:14} and
\eqref{eq:15}. When  $k\gamma=1$, on the other hand, one easily checks
that \eqref{eq:hg} implies $\int_{[\lambda^{-k}
  x,1]} \zeta_k(\mathrm du)
u^{\frac{1-k\varepsilon}{k}}=O(\lambda^{k\varepsilon}x^{-\varepsilon})$, entailing \eqref{eq:16}.
$\hfill \square$

\bigskip

\textbf{Proof of Lemma \ref{sec:proof-prop-refs}.} 
Since $\sum_{i\geq 1}s_i=1$ $\nu$-a.e., we have, 
\begin{align*}
\int_{\mathcal S} \nu(\mathrm d \mathbf s) \Bigg(\sum_{i\geq 1} g_k(xs_i^{-k})-g_k(x) \Bigg)^2 &\leq  2 \Bigg(\int_{\mathcal S} \nu(\mathrm d \mathbf s) \Bigg(\sum_{i\geq 2}( g_k(xs_i^{-k})-s_ig_k(x)) \Bigg)^2 \Bigg.  \\
&\qquad+ \int_{\mathcal S} \nu(\mathrm d \mathbf s) \left( g_k(xs_1^{-k})-s_1g_k(x) \right)^2 \mathbbm 1_{\{s_1\leq 1/2\}} \\
&\qquad+ \Bigg.\int_{\mathcal S} \nu(\mathrm d \mathbf s) \left( g_k(xs_1^{-k})-s_1g_k(x) \right)^2 \mathbbm 1_{\{s_1> 1/2\}} \Bigg). 
\end{align*}

{\bf  When $\gamma \in (0,1)$ and $k\gamma<1$,} by \eqref{eq:15}, the above upper bound is smaller than 
$$
 2 \kappa ^2 \cdot (x^{-\frac{1}{k}+\varepsilon})^2  \int_{\mathcal S} \nu(\mathrm d \mathbf s) \left(\left(\sum_{i\geq 2} s_i^{1-k \varepsilon} \right)^2+s_1^{2(1-k \varepsilon)} \mathbbm 1_{\{s_1 \leq 1/2\}}+s_1^2 (1-s_1)^{2\gamma} \mathbbm 1_{\{s_1> 1/2\}}  \right) 
$$
(we used that $\nu$-a.e. $s_i \leq 1/2$ for all $i\geq 2$) for some
finite $\kappa$ and $\varepsilon>0$. Taking
$\varepsilon$ smaller if necessary, so that  $2k \varepsilon <\eta$
(with the $\eta$ of assumption \eqref{eq:exp}), we claim that this
last integral is finite, entailing \eqref{eq:17}. To see this, first note
that
$$\int_{\mathcal S} \nu(\mathrm d \mathbf s) \left(s_1^{2(1-k
    \varepsilon)} \mathbbm 1_{\{s_1 \leq 1/2\}}+s_1^2 (1-s_1)^{2\gamma}
  \mathbbm 1_{\{s_1> 1/2\}}  \right)
\leq \nu(\{s_1\leq 1/2\})+\int_{\mathcal{S}}\nu(\mathrm d\mathbf{s})
(1-s_1)^{2\gamma}<\infty\, .$$
Indeed, $ \nu(\{s_1\leq 1/2\})\leq 2\int_{\mathcal S}(1-s_1) \nu(\mathrm
d\mathbf s)<\infty$, while \eqref{eq:hg} entails that $\int_{\mathcal S}(1-s_1)^{2\gamma} \nu(\mathrm
d\mathbf s)<\infty$. We finally control the last term in the integral
by the Cauchy-Schwarz inequality, using the fact that $\sum_{i\geq
  2}s_i\leq 1$, and assumption \eqref{eq:exp}: 
$$\int_{\mathcal S} \nu(\mathrm d\mathbf{s})\left(\sum_{i\geq 2}
  s_i^{1-k \varepsilon} \right)^2\leq \int_{\mathcal S} \nu(\mathrm d\mathbf{s})\sum_{i\geq 2}
  s_i^{1-2k \varepsilon} <\infty.$$

{\bf  When $\gamma=0$, }we proceed similarly with the bound given by \eqref{eq:14}, the only difference being that here we can use the following bound 
$$
\int_{\mathcal S} \nu(\mathrm d \mathbf s) \left( g_k(xs_1^{-k})-s_1g_k(x) \right)^2 \mathbbm 1_{\{s_1> 1/2\}}  \leq \kappa ^2 \cdot (x^{-\frac{1}{k}+\varepsilon})^2 \cdot \int_{\mathcal S} \nu(\mathrm d \mathbf s) s_1^{2-2k\varepsilon}, 
$$ 
where the integral is finite for $\varepsilon$ small enough since $\nu$ is finite and $s_1 \leq 1$. 

{\bf When $k\gamma=1$,} we proceed similarly, now with the
bound \eqref{eq:16}, which leads to 
\begin{multline*}
\int_{\mathcal S} \nu(\mathrm d \mathbf s) \Bigg(\sum_{i\geq 1} g_k(xs_i^{-k})-g_k(x) \Bigg)^2 \\
\leq 2\kappa^2 \cdot  x^{-\frac{2}{k}}  \int_{\mathcal S} \nu(\mathrm d \mathbf s)  \left(\Bigg(\sum_{i\geq 2} s_i |\log(s_i)| \Bigg)^2+s_1^2 | \log(s_1)|^{2}  \mathbbm 1_{\{s_1 \leq 1/2\}} +s_1^2 (1-s_1)^{2\gamma}  \mathbbm 1_{\{s_1> 1/2\}} \right) 
\end{multline*}
and again the integral is finite under \eqref{eq:hg} and
\eqref{eq:exp}. 
$\hfill \square$

\section{Examples}\label{sec:examples}

\subsection{Cases with finite $\nu$: Dirichlet
  fragmentations}\label{sec:cases-with-finite} 

When $\nu$ is finite, \eqref{eq:hg} with $\gamma=0$ is automatically satisfied
with $c_\nu=\nu(\mathcal{S})$ (which we will normalise to $1$
in the examples below). Let us discuss some examples of this
situation. 

Fix some $m\geq 2$ and a family $(a_1,\ldots,a_m)$ of positive
numbers. The Dirichlet distribution $\mathrm{Dir}(a_1,\ldots,a_m)$ is
the probability distribution on the simplex $S_{m-1}=\{(x_1,\ldots,x_m)\in
\R_+^m:x_1+\cdots+x_m=1\}$
with density $B(a_1,\ldots,a_m)\prod_{i=1}^mx_i^{a_i-1}$ with respect
to the uniform probability measure $\Delta_m$, where
$B(a_1,\ldots,a_m)=\Gamma(a_1)\cdots \Gamma(a_m)/\Gamma(a_1+\cdots+a_m)$. Let us now consider the measure
$\nu_{(a_1,\ldots,a_m)}$ that is the push-forward of the Dirichlet distribution
$\mathrm{Dir}(a_1,\ldots,a_m)$ by  the mapping
$$\mathbf{x}=(x_1,\ldots,x_m)\mapsto (x_{(1)},\ldots,x_{(m)},0,\ldots)\, ,$$
where $x_{(1)}\geq x_{(2)}\geq \ldots \geq x_{(m)}$ is the
nonincreasing rearrangement of $x_1,\ldots,x_m$. For instance, the
case $m=2$ and $a_1=a_2=1$ corresponds to  successively splitting
intervals in two subintervals at uniformly 
random locations.

In this model, \eqref{eq:exp} and  \eqref{eq:gamma0} hold trivially, and we are in the subcritical case of Theorem \ref{sec:thm-sub-gammapos} for all $k\geq 2$. The constant $C^{\mathrm{sub}}_{\nu_{(a_1,\ldots,a_m)}}(k)$ is not particularly nice, but reasonably
explicit: if $\Delta_m$ is the uniform probability measure on the
simplex of dimension $m-1$, then 
\begin{equation}
  \label{eq:33}
  C_{\nu_{(a_1,\ldots,a_m)}}^{\mathrm{sub}}(k)=\frac{\int_{S_{m-1}}
  \mathrm{Dir}(a_1,\ldots,a_m)(\mathrm d\mathbf{x})
  (1-\sum_{i=1}^mx_i^k)^{1/k}}{\frac{\Gamma'(a_1+\cdots+a_m+1)}{\Gamma(a_1+\cdots+a_m+1)}-\frac{1}{a_1+\cdots+a_m}\sum_{i=1}^m\frac{\Gamma'(a_i+1)}{\Gamma(a_i)}}\, .
\end{equation}

\medskip

{\remark It would be interesting to study the case of the $m$-ary
  fragmentation measure $\nu(s_i\in \mathrm dx)=\delta_{1/m}(\mathrm
  dx)$ for $1\leq i\leq m$, and $\nu(s_i\in \mathrm
  dx)=\delta_0(\mathrm dx)$ for $i>m$. In this situation, the potential measure
  $U$ is lattice, and our results do not apply. 
  In this apparently very simple
  case, the urn scheme is completely explicit, but $S_x$ has an
  oscillatory behavior. The strong law of \cite[Section 5]{karlin67}
  apply and show that $N^{(\nu)}_n(k)/\mathbb{E}[N^{(\nu)}_n(k)]$
  converge a.s.\ in this case, but the behavior of
  $\mathbb{E}[N^{(\nu)}_n(k)]$ can be quite complicated. } 

\subsection{Cases with infinite $\nu$}

\subsubsection{Stable trees}

Indexed by a parameter $\beta \in (1,2]$, the stable trees introduced in
\cite{LGLJ98,DuLG02} generalize the Brownian CRT to heavy tailed settings, with an important role in branching and random graphs theories. The stable tree of exponent 
$\beta=2$ is simply a version of the
Brownian CRT, and we use here the convention that it is a version of the Brownian CRT considered in Section \ref{sec:brownian-tree-case} where the distances are multiplied by $2^{1/2}$, that is, with a dislocation measure $\nu_2$ equals to $2^{-1/2}$ times the measure defined in (\ref{nu_brownien}).  In \cite{mierfmoins}, it was proved that the stable tree of exponent $\beta \in (1,2)$ is
also a fragmentation tree, now with index of self-similarity
$\beta^{-1}-1$ and dislocation measure $\nu_{\beta}$ given by
\begin{equation}
\label{nu_stable}
 \int_{\mathcal{S}^\downarrow}\nu_\beta(\mathrm{d}\mathbf{s})f(\mathbf{s})=\frac{ \Gamma(1-1/\beta)}{\Gamma(-\beta)}\cdot\mathbb{E}\left[T^{(1/\beta)}_1f\left(\frac{\Delta^{\downarrow}
      T^{(1/\beta)}_{[0,1]}}{T_1^{(1/\beta)}}\right)\right]\, ,
      \end{equation}
where $(T^{(1/\beta)}_x,x\geq 0)$ is a stable subordinator of Laplace exponent
$\lambda^{1/\beta}$, and $\Delta^{\downarrow}
      T^{(1/\beta)}_{[0,1]}$ is the sequence of its jumps over the
      interval $[0,1]$, ranked by decreasing order of magnitude.

We can treat in the same go the multiple of the Brownian CRT  and the
stable trees with exponent $\beta\in (1,2)$. Indeed, we know from
\cite[p.440]{mierfmoins} 
that the Laplace exponent $\phi_{\beta}$ of the associated subordinator (see Section \ref{sec:tagged-fragment}) of the $\beta$-model is
$$\phi_{\beta}(q)=\frac{\beta \Gamma(q+1-1/\beta)}{\Gamma(q)}\, .$$
In
particular,  $\phi_{\beta}(q)\sim \beta \, q^{1-1/\beta}$ as $q\rightarrow
\infty$, using Stirling's formula, hence \eqref{eq:hg} holds with
$\gamma=1-1/\beta$ and
$c_{\nu_\beta}=\beta/\Gamma(1/\beta)$. Moreover,
$\phi_{\beta}'(0+)=\beta \Gamma(1-1/\beta)$. We also note, as in
\cite[p.4345]{qshi15}, that it yields the explicit expression 
$(\beta \Gamma(1-1/\beta))^{-1}\cdot (1-e^{-y})^{-1/\beta}) 
\mathrm dy$ for the associated potential measure, but we will not need
this.

For $\beta=2$, the dislocation measure is binary, so that \eqref{eq:exp} holds
automatically. But in fact, rewriting $\phi_{\beta}$ as
$$\phi_{\beta}(q)= \frac{\beta q\Gamma(q+1-1/\beta)}{\Gamma(q+1)}\, ,$$
we see that $\phi_{\beta}$ can be analytically continuated in a neighborhood
of $0$. By the discussion around \eqref{eq:27}, this shows that
\eqref{eq:exp} holds for all $\beta \in (1,2]$.

In this setting, Theorem \ref{sec:thm-super} and Theorem \ref{sec:thm-sub-gammapos} read as follows (we give the statement for the number of ancestors $N^{(\nu_{\beta})}_n(k)$, the statement for $N^{(\nu_{\beta})}_{n,r}(k)$ is easily adapted). We slightly change perspective by fixing the integer $k$ and letting $\beta$ varies.

\vspace{0.15cm}

\begin{proposition}
Fix an integer $k\geq 2$. Then almost surely and in $L^2$, as $n \rightarrow \infty$,
\vspace{-\topsep}
\begin{itemize}
\item if $\beta>k/(k-1)$,  then $n^{-1+1/\beta}N^{(\nu_{\beta})}_n(k) \to
  \beta k^{1-1/\beta} \cdot \mathcal{A}^{(\nu_\beta)}_{k},~$ where
  $\mathcal{A}_{k}^{(\nu_\beta)}$ is the area of a
  $(1-k(1-1/\beta),\nu_\beta)$-fragmentation tree 
  \emph{(supercritical case)}, 
\item if $\beta=k/(k-1)$, then
  $(n^{1/k}\log(n))^{-1}N^{(\nu_{\beta})}_n(k)   \to
  k^{\frac{1}{k}-1} /\Gamma(1/k)~$  \emph{\ (critical case)},
\item if $\beta<k/(k-1)$, then $$n^{-1/k}N^{(\nu_{\beta})}_n(k) \to  \frac{\Gamma(1-1/k)}{|\Gamma(1-\beta)|} \cdot  
 \mathbb{E}\left[T^{(1/\beta)}_1\left( 1-\sum_{i=1}^{\infty}\left(
       \frac{\Delta^{(1/\beta)}_{i}}{T_1^{(1/\beta)}}\right)^k\right)^{1/k}\right]\, ,$$
   where $T^{(1/\beta)}$ is the $1/\beta$-stable subordinator introduced in (\ref{nu_stable}) \emph{(subcritical case)}.
\end{itemize}

\end{proposition}

Observe that in the subcritical case, whatever the value of $k\geq 2$, $\beta$ cannot be equal to 2.

\subsubsection{Ford's trees}

A planted binary tree is a rooted tree in which all vertices have
degree $1$ or $3$, and the root vertex has degree $1$. An edge in a
binary tree is called \emph{external} if it is incident to a vertex of
degree $1$ that is distinct from the root vertex, and is called
\emph{internal} otherwise. Note that a planted binary trees with
$n\geq 2$
external edges must have $n-1$ internal edges. 
Ford's model of growing trees is a Markov
chain $(T_n,n\geq 2)$ on the set of planted binary trees, depending on
a parameter $a\in (0,1)$, and defined as
follows. We let $T_2$ be planted binary tree with two external edges, and one internal
edge. At step $n\geq 2$, an edge of $T_n$ is selected at random, with
probability proportional to $a$ if the edge is internal, and with
probability proportional to $1-a$ if the edge is external. We then
graft a new external edge to the middle of the selected edge. More
formally, we subsitute to the selected edge, say $\{x,y\}$, where $x$
and $y$ are vertices of $T_n$, a star-graph
$(\{x,y,x',y'\},(\{x,x'\},\{x',y\},\{x',y'\})$, where $x',y'$ are two
new vertices, not in $T_n$. We call $T_{n+1}$ the resulting tree,
which obviously has $n+1$ external edges. Note that for $a=1/2$, the
above Markov chain is known as Rémy's algorithm, and generates at time
$n$ a uniformly random binary tree with $n$ exterior edges (when
labeled in order of appearance).  

It was shown in \cite[Section 5.2]{HaMiPiWi08} that some versions of the trees $T_n,n\geq 2$ can in fact be recovered by
a simple sampling procedure of a self-similar fragmentation
tree. Namely, letting $\nu_a$ be the measure on $\mathcal{S}$ such
that $\nu_a(\{\mathbf{s}\in \mathcal{S}:s_1+s_2<1\})=0$ and 
$$\frac{\nu_a(s_1\in \mathrm dx)}{\mathrm
  dx}=\frac{a(\Gamma(1-a))^{-1}}{(x(1-x))^{a+1}}+\frac{2(1-2a)(\Gamma(1-a))^{-1}}{(x(1-x))^a}\, ,\quad
x\in (1/2,1)\, ,$$
we can let $(\mathcal{T}_a,d_a,\rho_a,\mu_a)$ be the self-similar
fragmentation tree with dislocation measure $\nu_a$, and
self-similarity index $-a$. Then, if $x_1,x_2,\ldots$ is an i.i.d.\
sample of points distributed according to $\mu_a$ and $T'_n$ is the
combinatorial skeleton of the subtree of $\mathcal{T}_a$ spanned by
the root and the points $\{x_1,\ldots,x_n\}$, it holds that $T'_n$ has
same distribution as $T_n$, for each $n$, even though 
the distributions of the sequences $(T'_n,n\geq 2)$ and $(T_n,n\geq 2)$ are not equal. Consequently, if we group, for example, the leaves of the tree two by two (leaving one unpaired if $n$ is odd), the number of different most recent common ancestors is distributed as $N_{\lfloor n/2 \rfloor}^{(\nu_a)}(2)$. The proposition below therefore provides the behaviour in distribution of the  ancestor-counting variables related to the tree $T_n$. 

The binary measure $\nu_a$ clearly satisfies (\ref{eq:exp}) and (\ref{eq:gamma0}). Moreover one checks that 
$$
\phi_a(q)= \frac{\Gamma(q+1-a)\Gamma(q+2)}{\Gamma(q) \Gamma(q+3-2a)}.
$$
In particular $\phi_a(q)\underset{q \rightarrow \infty}\sim q^a$ by Stirling's formula, and therefore (\ref{eq:hg}) holds with $\gamma=a$ and $c_{\nu_a}=1/\Gamma(1-a)$. Last, $\phi_a'(0+)=\Gamma(1-a)/\Gamma(3-2a)$. In this case, our results therefore resume as follows:

\begin{proposition}
  \label{sec:fords-trees}
Let $k\geq 2$ be a fixed integer.  Then almost surely and in $L^2$, as $n \rightarrow \infty$,
\vspace{-\topsep}
   \begin{itemize}
  \item if $a>1/k$, then
    $n^{-a}N^{(\nu_a)}_n(k)\to k^a\mathcal{A}_k^{(\nu_a)}$, where
    $\mathcal{A}_k^{(\nu_a)}$ is the area of a
    $(1-ka,\nu_a)$-fragmentation tree  \emph{(supercritical case)}, 
   \item if $a=1/k$, then
     $(n^{1/k}\log(n))^{-1}N^{(\nu_a)}_n(k)\to k^{\frac{1}{k}-1} \Gamma(3-2/k)/\Gamma(1-1/k)$ \emph{(critical case)}, 
   \item
     if $a <1/k$, then $n^{-1/k}N^{(\nu_a)}_n(k)\to\Gamma(1-1/k) \cdot
    C^{\mathrm{sub}}_{\nu_a}$ \emph{(subcritical case)}.
  \end{itemize}
\end{proposition}

\smallskip

{\remark Although we do not need it for our purposes, we identified explicitly the potential measure associated with Ford's model while working on this problem. We give it here since it may have its own interest.}

\begin{proposition} Let $a \in (0,1)$. The potential measure associated to a fragmentation tree with dislocation measure $\nu_a$  through the relation (\ref{eq:32}) (with $\lambda=0$) is absolutely continuous with a density defined by
$$f_a(t)=g_a(e^{-t}), \quad t>0,$$
where
\begin{eqnarray*}
g_a(x)&=& \frac{1}{\Gamma(a)} x^{3-2a}(1-x)^{a-1} \sum_{n=0}^{\infty} \frac{(2)_n(1-a)_n}{(a)_n} \cdot  \frac{(1-x)^n}{n!}
\end{eqnarray*}
where $(u)_n=u(u+1)\ldots(u+n-1)$ is the Pochhammer symbol. 
\end{proposition}

\textbf{Proof.} We use Gauss's summation theorem:
$$
\sum_{n=0}^{\infty} \frac{(x)_n(y)_n}{(z)_n} \cdot \frac{1}{n!} =\frac{\Gamma(z) \Gamma(z-x-y)}{\Gamma(z-x)\Gamma(z-y)} \quad \text{when } z>x+y,
$$
together with Fubini-Tonelli's theorem to see that the function
$$
q \in (0,\infty) \mapsto \frac{1}{\phi(q)}=\frac{\Gamma(q) \Gamma(q+3-2a)}{\Gamma(q+1-a)\Gamma(q+2)}
$$
is the Mellin transform of $g_a$, and therefore the Laplace transform of $f_a$, as required.
$\hfill \square$

\subsection{Infinite Beta-type dislocation measures}

We consider extensions to infinite dislocation measures of the model of Dirichlet fragmentations of Section \ref{sec:cases-with-finite}. To simplify, we focus on dislocations into $m=2$ pieces. Let $a>-1,b>-1$ be two parameters and consider the \emph{binary} dislocation measure characterized by the distribution of its largest fragment as follows:
\begin{equation}
\label{def:nubeta}
\nu_{(a,b)}(s_1 \in \mathrm dx)=\left(x^{a-1}(1-x)^{b-1}+x^{b-1}(1-x)^{a-1} \right)\mathbbm1_{\{1/2 \leq x <1\}} \mathrm dx.
\end{equation}
By symmetry, we may assume that $a \geq b >-1$. Note that  $\int_{\mathcal S} (1-s_1) \nu_{(a,b)}(\mathrm d\mathbf s)$ is indeed finite and that the measure $\nu_{(a,b)}$ is itself finite if and only if $b>0$, resuming then to the situation of Section \ref{sec:cases-with-finite} with dislocations into two pieces according to a $\mathrm{Beta}(a,b)$ distribution. 

This extension to infinite Beta-type dislocation measures encompasses
the scaling limits of Aldous' $\beta$-splitting trees, when
$a=b=\beta+1$ for some $\beta  \in (-2,\infty)$. Aldous
$\beta$-splitting trees have been introduced in  \cite{aldous96clad}
as theoretical models for phylogenetic trees, see
\cite{lambert17,steel2016} for overviews on that topic. The  $\beta$-splitting trees are discrete rooted trees with $n$ leaves coding the evolution of  ``clades'', where clades are recursively split into sub-clades, with the rule that a clade of  $k$ leaves is split into sub-clades containing  $i$ and  $k-i$ leaves at a rate proportional to  $\frac{\Gamma(\beta+i+1)\Gamma(\beta+k-i+1)}{\Gamma(i+1)\Gamma(n-i+1)}$. When $\beta \in (-2,-1)$, the height of the tree is then proportional to $n^{-\beta-1}$ and the limit of the rescaled tree has been identified in \cite{HaMiPiWi08} has a fragmentation tree with parameters $(\beta+1,\nu_{(\beta+1,\beta+1)})$. In particular, when $\beta=-3/2$, one recovers the Brownian tree up to a multiplicative constant. For any $\beta \in (-2,-1)$, by considering an infinite sample of i.i.d. leaves of the fragmentation tree, and, for each $n$,  the combinatorial skeleton of the subtree spanned  by the root and the $n$ first sampled leaves, one recover a version of the $\beta$-splitting tree with $n$ leaves. Proposition \ref{prop:betasplitting} below therefore concerns the  ancestor-counting variables for both the discrete and continuous models.
When $\beta > -1$, the dislocation measure becomes finite, which
simplifies a lot the structure of the genealogy. The critical case
$\beta=-1$ is of notable interest and was recently studied in
\cite{aldousjanson2024}.

Back to the general model, we note that for any  $a \geq b >-1$, the binary measure $\nu_{(a,b)}$ defined by (\ref{def:nubeta}) satisfies the assumptions (\ref{eq:exp}) and (\ref{eq:gamma0}). Moreover, \newline
$\begin{array}{ll} \bullet \text{ when }b>0, & \nu_{(a,b)}(\mathcal S)= \frac{\Gamma(a)\Gamma(b)}{\Gamma(a+b)}\\ \bullet  \text{ when }b=0, & \left\{\begin{array}{ll} \nu_{(a,b)}(s_1<1-x) \underset{x \downarrow 0}{\sim} |\log (x)| & \text{ if }a >0 \\ \nu_{(a,b)}(s_1<1-x) \underset{x \downarrow 0}{\sim} 2|\log (x)|  & \text{ if }a =0 \end{array}\right.\\ \bullet \text{ when }b<0 & \left\{\begin{array}{ll} \nu_{(a,b)}(s_1<1-x) \underset{x \downarrow 0}{\sim}  |b|^{-1}x^{b} & \text{ if }a >b \\  \nu_{(a,b)}(s_1<1-x) \underset{x \downarrow 0}{\sim} 2|b|^{-1}x^{b} & \text{ if } a=b . \end{array}\right. \end{array}$

Therefore, our assumption (\ref{eq:hg}) holds with
  $\gamma=\max(-b,0)$, except when $b=0$. This latter case 
  corresponds to an extension of (\ref{eq:hg}) to a regular variation
  situation.  Although we believe that our results could be extended
  in general to a regularly varying version of (\ref{eq:hg}), we have
  not checked it properly. However, for the present model when $b=0$,
  corresponding to a subcritical regime, we did check that all steps
  of our proof are indeed valid. To summarise, one can check using
  \cite[Theorem 27.7]{sato13}  that the
  tagged fragment subordinator has an absolute continuous density for
  every positive time, which warrants the use of the renewal theory and the concentration results of Section \ref{sec:concentration1}. The main differences lie in the proof of Lemma \ref{lem:increments_g}, where typically $U([0,x])$ is now proportional to $|\log(x)|$ when $x \downarrow 0$, instead of a constant for the usual (\ref{eq:hg}) assumption when $\gamma=0$. 

Besides, one easily checks that 
\begin{eqnarray*}
\phi'_{(a,b)}(0+)&=&\int_0^1 \left( |\log(u)| u+ |\log(1-u)| (1-u) \right) u^{a-1}(1-u)^{b-1} \mathrm du \\
&=& \frac{\Gamma(a)\Gamma(b)}{\Gamma(a+b+1)}
    \left(\frac{\Gamma'(a+b+1)}{\Gamma(a+b)}-\frac{\Gamma'(a+1)}{\Gamma(a)}-\frac{\Gamma'(b+1)}{\Gamma(b)}
    \right)\, .
\end{eqnarray*}
Fix an integer $k\geq 2$. Then, setting
$$
C^{\mathrm{sub}}_{\nu_{(a,b)}}(k)=\frac{\int_0^1 \left(1-u^k-(1-u)^k\right)^{1/k}u^{a-1}(1-u)^{b-1} \mathrm du}{\phi_{(a,b)}'(0+)}
$$
(which is \eqref{eq:33} with $m=2$ and $(a_1,a_2)=(a,b)$ when $b>0$,
as it should) and
$$
C^{\mathrm{cr}}_{\nu_{(a,b)}}(k)= \frac{(1+\mathbbm 1_{\{a=b\}})k^{\frac{1}{k}}}{\phi_{(a,b)}'(0+)},
$$
our results read on this model as follows:

\begin{proposition} 
\label{prop:betasplitting}
Almost surely and in $L^2$, as $n\to\infty$,
\vspace{-\topsep}
\begin{itemize}
  \item if $b<-1/k$,  then $n^b N^{\nu_{(a,b)}}_n(k) \to
  |\Gamma(b)|(1+\mathbbm 1_{\{a=b\}}) k^{|b|}  \cdot
  \mathcal{A}^{(\nu_{(a,b)})}_{k} $, where
  $\mathcal{A}^{(\nu_{(a,b)})}_{k}$ is \linebreak the area of a
  $(1+kb,\nu_{(a,b)})$-fragmentation tree \emph{(supercritical case)},
  \item if $b=-1/k$, then $(n^{1/k}\log(n))^{-1}N^{\nu_{(a,b)}}_n(k) \to
  \Gamma(1-\frac{1}{k}) \cdot C^{\mathrm{cr}}_{\nu_{(a,b)}}(k)$ \emph{(critical case)},
\item if $b>-1/k$, then $n^{-1/k}N^{\nu_{(a,b)}}_n(k) \to
  \Gamma(1-\frac{1}{k}) \cdot C^{\mathrm{sub}}_{\nu_{(a,b)}}(k)$
  \emph{(subcritical case)}. 
\end{itemize}
\end{proposition}

\bibliographystyle{siam}
\bibliography{Bib}

\end{document}